\documentclass[10pt]{amsart}

\usepackage{hyperref}

\usepackage[dvips]{graphicx}
\usepackage[latin1]{inputenc}

\input xy
\xyoption{all}

\usepackage{amsmath}
\usepackage{amsfonts}
\usepackage{amssymb}
\usepackage{amsthm}
\usepackage{indentfirst}

\usepackage[T1]{fontenc}

\newcommand{\comment}[1]{}

\DeclareMathOperator\Perm{Perm}

\DeclareMathOperator\Decr{\mathcal{D}}

\newtheorem{Def}{Definition}[subsection]
\newtheorem{Defprop}[Def]{Definition-Proposition}
\newtheorem{Th}{Theorem}[subsection]
\newtheorem{Defth}[Th]{Definition-Theorem}
\newtheorem{Prop}[Th]{Proposition}
\newtheorem{Lem}[Th]{Lemma}
\newtheorem{Corol}[Th]{Corollary}
\newtheorem{Conj}{Conjecture}
\newtheorem{Ques}[Conj]{Question}\theoremstyle{remark}
\newtheorem{Rem}{Remark}
\newtheorem{example}{Example}

\comment{%version française
\usepackage[francais]{babel}
\newtheorem{Def}{Définition}
\newtheorem{Defprop}[Def]{Définition-Proposition}
\newtheorem{Th}{Théorème}
\newtheorem{Defth}[Th]{Définition-Théorème}
\newtheorem{Prop}[Th]{Proposition}
\newtheorem{Lem}[Th]{Lemme}
\newtheorem{Corol}[Th]{Corollaire}

\theoremstyle{remark}
\newtheorem{Rem}{Remarque}
}

\newcommand{\C}{\mathbb{C}}

\newcommand{\Z}{\mathbb{Z}}
\newcommand{\N}{\mathbb{N}}
\newcommand{\M}{\mathcal{M}}
\newcommand{\ck}{(1 2 \ldots k)}
\newcommand{\s}[1]{\mathbf{#1}}
\newcommand{\dec}[1]{#1,\overline{#1}}

\title[Positivity of Kerov's polynomials]{Combinatorial interpretation and positivity of Kerov's character polynomials}
\author{Valentin F\'eray}
\address{The Gaspard--Monge Institut of electronique and computer science,
University of Marne-La-Vall\'ee Paris-Est,
77454 Marne-la-Vall\'ee Cedex 2, France}
\email{feray@univ-mlv.fr}


\begin{document}

\begin{abstract}
Kerov's polynomials give irreducible character values in term of the free cumulants of the associated Young diagram. We prove in this article a positivity result on their coefficients, which extends a conjecture of S. Kerov.\\

Our method, through decomposition of maps, gives a description of the coefficients of the $k$-th Kerov's polynomials using permutations in $S(k)$. We also obtain explicit formulas or combinatorial interpretations for some coefficients. In particular, we are able to compute the subdominant term for character values on any fixed permutation (it was known for cycles).
\end{abstract}

\maketitle

\section{Introduction}
\subsection{Background}\label{parbackgr}

\subsubsection{Representations of the symmetric group}
Representations theory of the symmetric group $S(n)$ is a very ancient research field in mathematics. Irreducible representations of $S(n)$ are
indexed by partitions\footnote{Non-increasing sequences of non-negative
  integers of sum $n$.} $\lambda$ of $n$, or equivalently by Young
diagrams of size $n$. The associated character can be computed thanks
to a combinatorial algorithm, but unfortunately it becomes quickly
combersome when the size of the diagram is large and does not help to
study asymptotic behaviours.

\subsubsection{Free cumulants}
To solve asymptotic problems in representation theory of the symmetric
groups, P. Biane introduced in \cite{Bi1} the free cumulants
$R_i(\lambda)$ (of the transition measure) of a Young
diagram\footnote{The transition measure of a Young diagram is a
  measure on the real line introduced by S. Kerov in \cite{Ke}. Its free
  cumulants are a sequence of real numbers associated to this
  measure. The denomination comes from free probability theory, see
  \cite{Bi1} for more details.}. Asymptotically, the character value
and the classical operation on representations can be easily described
with with free cumulants: 
\begin{itemize}
\item Up to a good normalisation, the $l+1$-th free cumulant is the leading term of the character value on the cycle $(1 \ldots l)$.
\item Typical large Young diagrams (according to the Plancherel distribution) have, after rescaling, all their free cumulants, excepted from the second one, very close to zero.
\item Almost all the diagrams appearing in an elementary operation on
  irreducible representations (like restriction, tensor product) have
  free cumulants very close to specific values, which can be easily
  computed from the free cumulants of the original diagram(s).
\end{itemize}
So the free cumulants form a good way to encode the informations contained in a Young diagram.\\

\subsubsection{Kerov's polynomials}
It is natural to wonder if there are exact expressions of character value in terms of free cumulants. Kerov's polynomials give a positive answer to this question for character values on cycles (they appear first in a paper of P. Biane \cite[Theorem 1.1]{Bi2} in 2003).
Unfortunately, their coefficients remain very mysterious. A lot of
work has been done to understand them
(\cite{Bi2},\cite{\'Sn2},\cite{GR},\cite{Bi4},\cite{R\'S},\cite{La}):  a general,
but exploding in complexity, explicit formula and a combinatorial interpretation for linear terms in free cumulants have been found.\\

The positivity of the coefficients of Kerov's polynomial has been observed by numerical computations (\cite{Bi2},\cite{GR}) and was conjectured by S. Kerov. The main result of this paper is a positive answer to this conjecture.\\

\subsubsection{Multirectangular Young diagrams}
\begin{figure}
\begin{center}\includegraphics{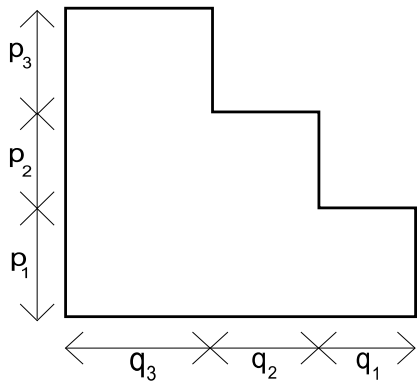}\caption{\footnotesize Young diagram associated to sequences $\s{p}$ and $\s{q}$ (french convention)}\label{figpart}
\end{center}
\end{figure}

We use in this paper a new way to look at Young diagrams, initiated by R. Stanley in \cite{St1}. In this paper, he proved a 
nice combinatorial formula for character values, but only for Young diagrams of rectangular shape. To generalize it,
we have to look at any Young diagram as a superposition of rectangles as in figure \ref{figpart}. With this description, Stanley's
formula has been recently generalized (see \cite{St2},\cite{F\'e}).\\

The complexity of this general formula depends only on the size of the
support of the permutation (and not of the size of the permutation itself!). As
remarked in \cite{F\'S}, it is useful to reformulate it with the
notion of bipartite graph associated to a pair of permutations. This
bipartite graph has in fact a canonical map structure\footnote{For
  some pairs of permutations, this structure was introduced by
  I.P. Goulden and D.M. Jackson in \cite{GJ}.}, which is central
here.\\

In this paper, we link these two recent developments. This gives a new combinatorial interpretation of the coefficients, proving Kerov's conjecture.

\subsection{Normalized character}
If $\sigma$ is a permutation in $S(k)$, let $C(\sigma)$ be the
partition of the set $[k]:=\{1,\ldots,k\}$ in orbits under the action of
$\sigma$. The type of $\sigma$ is, by definition, the partition $\mu$
of the integer $k$ whose parts are the length of the cycle of
$\sigma$. The conjugacy classes of $S(k)$ are exactly the sets of partition of a given type.\\

By definition, for $\mu \vdash k$ and $\lambda \vdash n$ with $k \leq n$, the normalized cha\-racter value is given by equation:
\begin{equation}
\Sigma_{\mu}(\lambda):= \frac{n(n-1)\ldots (n-k+1) \chi^\lambda(\sigma)}{\chi^\lambda(Id_n)},
\end{equation}
where $\sigma$ is a permutation in $S(k)$ of type $\mu$ and $\chi^\lambda$ is the character value of the irreducible representation associated to $\lambda$ (see \cite{McDo}). Note that we have to identify $\sigma$ with its image by the natural embedding of $S(k)$ in $S(n)$ to compute $\chi^\lambda(\sigma)$.\\

\comment{
One of their nice properties is that they give the leading term of the character value on cycles for balanced Young diagrams ($\sum p_i+\sum q_i = O(\sqrt{|\lambda|})$), more precisely (see P. Biane, \cite[Theorem 1.3]{Bi1}):
\begin{equation}\label{careqcum}
\Sigma_k(\lambda)=R_{k+1}(\lambda)+ O\big(|\lambda|^{\frac{k-1}{2}}\big)
\end{equation}
}

\subsection{Minimal factorizations and non-crossing partitions}\label{subsectncp}
Non-crossing partitions and in particular, their link with minimal factorizations of a cycle, are central in this work. This paragraph is devoted to definition and known results in this domain. For more details, see P. Biane's paper \cite{Bncp}.

\begin{Def}
A crossing of a partition $\pi$ of the set $[j]$ is a quadruple $(a,b,c,d) \in [j]^4$ with $a < b < c < d$ such that
\begin{itemize}
\item $a$ and $c$ are in the same part of $\pi$;
\item $b$ and $d$ are in the same part of $\pi$, different from the one containing $a$ and $c$.
\end{itemize}
A partition without crossings is called a non-crossing partition. The
set of non-crossing partitions of $[j]$ is denoted NC(j) and can be
endowed with a partial order structure (by definition, $\pi \leq \pi'$ if every part
of $\pi$ is included in some part of $\pi'$). 
\end{Def}

The partially ordered set (poset) $NC(j)$ appears in many domains: we will use its connection with the symmetric group.\\

Let us consider the following length on the symmetric group $S(j)$:
denote by $l(\sigma)$ the minimal number $h$ of transpositions needed
to write $\sigma$ as a product of transpositions $\sigma = \tau_1
\ldots \tau_h$. One has:
\begin{eqnarray*}
l(Id_j) & = & 0,\\
l(\sigma^{-1}) & = & l(\sigma),\\
l(\sigma \cdot \sigma') & \leq & l(\sigma) + l(\sigma').
\end{eqnarray*}
We consider the associated partial order on $S(j)$: by definition,
$\sigma \leq \sigma'$ if $l(\sigma') = l(\sigma) + l(\sigma^{-1}
\sigma')$. It is easy to prove that 
\begin{itemize}
\item $Id_j$ is the smallest element ;
\item for any $\sigma$, one has $l(\sigma)=j-|C(\sigma)|$.
\end{itemize}
So, if we denote by $(1 \ldots j)$ the cycle sending 1 onto
2, 2 onto 3, etc\ldots, one has
$$\sigma \leq (1 \ldots j) \Longleftrightarrow |C(\sigma)| + |C(\sigma^{-1} (1 \ldots j))|=j+1.$$

If $\sigma \leq \sigma'$, let us consider the interval $[\sigma;\sigma']$ which is by definition the set $\big\{ \tau \in S(k) \text{ s.t. }\sigma \leq \tau \leq \sigma' \big\}$. In his paper \cite[section 1.3]{Bncp}, P. Biane gives a combinatorial description of these intervals:\\

\begin{Prop}[Isomorphism with minimal factorizations]\label{bijpermncp}
The map
\begin{eqnarray*}
[Id_j ; (1 \ldots j)] & \longrightarrow & NC(j) \\
\sigma & \mapsto & C(\sigma)
\end{eqnarray*}
is a poset isomorphism.
\end{Prop}

Here is the inverse bijection: to a non-crossing partition $\tau$ of
$[j]$, we associate the permutation $\sigma_\pi \in S(j)$, where $\sigma_\pi(i)$ is the next element in the same part of $\pi$ as $i$ for the cyclic order $(1,2,\ldots,j)$.\\

Since the order is invariant by conjugacy, every interval $[Id_j ;
  c]$, where $c$ is a full cycle, is isomorphic as poset to a
non-crossing partition set. More generally, if $\sigma$ is a permutation in $S(j)$,
$$[Id_j;\sigma] \simeq \prod_{i=1}^{|C(\sigma)|} NC(j_i),$$ where the $j_i$'s are the number of elements of the cycles of $\sigma$. This result gives a description of all intervals of the symmetric group since, if $\sigma \leq \sigma'$, we have $[\sigma;\sigma'] \simeq [Id;\sigma^{-1} \sigma']$.

\subsection{Kerov's polynomials}\label{subsectkerpol}
We look for an expression of the normalized cha\-racter value in terms
of free cumulants. In the case where $\mu$ has only one part
($\mu=(k), \sigma =(1 \ldots k)$), P. Biane shows\footnote{P. Biane attributes this result to S. Kerov.} in \cite{Bi2} that:\\

\begin{Defth} For any $k \geq 1$, there exists a polynomial $K_k$, called $k-$th Kerov's polynomial, with integer coefficients, such that, for every Young diagram $\lambda$ of size bigger than $k$, one has:
\begin{equation}\label{Kerov}
\Sigma_k(\lambda)=K_k(R_2\big(\lambda), \ldots, R_{k+1}(\lambda)\big).
\end{equation}
\end{Defth}

\vspace{-.7 cm}

\begin{minipage}[t]{.4\linewidth}
   \begin{eqnarray*}
   \text{Examples:}\\
   \Sigma_1 & = & R_2 ;\\
   \Sigma_2 & = & R_3 ;
   \end{eqnarray*}
\end{minipage}
\begin{minipage}[t]{.6\linewidth}
   \begin{eqnarray*}
   \Sigma_3 & = & R_4 + R_2 ;\\
   \Sigma_4 & = & R_5 + 3R_3 ;\\
   \Sigma_5 & = & R_6 + 15 R_4 + 5 R_2^2 + 8 R_2.
   \end{eqnarray*}
\end{minipage}\\

\vspace{0,2cm}

Our main result is the positivity of the coefficients of Kerov's polynomials. This result was conjectured by S. Kerov (according to P. Biane, see \cite{Bi2}).\\

\begin{Th}[Kerov's conjecture]\label{mainth}
For any integer $k \geq 1$, the polynomial $K_k$ has non-negative coefficients.
\end{Th}

Our proof gives a (complicated) combinatorial interpretation of the coefficients and allows us to compute some of them.\\

\subsubsection{High graded degree terms}

\begin{Th}\label{thgenre1}
Let $j_1,\ldots,j_t$ be non negative integers such that $\sum\limits_i j_i = k-1$. The coefficient of $\prod\limits_i R_{j_i}$ in $K_k$ is
\begin{equation}\label{eqgenre1}
\frac{(k-1)k(k+1)}{24} |\Perm(\s{j})| \prod\limits_i (j_i-1),
\end{equation}
where $\Perm(\s{j})$ is the set of sequences equal to $\s{j}$ up to a permutation ($|\Perm(\s{j})|=\frac{t!}{m_2! \ldots m_{k-1}!}$ is the multinomial coefficient of the $m_l$'s, where $m_l$ is the number of $j_i$ equal to $l$). 
\end{Th}
This theorem gives an explicit formula for the term \comment{$\Sigma{k;2}$} of graded degree $k-1$ in $K_k$, which is the subdominant term for character values on a cycle. It has already been proved in two different ways by I.P. Goulden and A. Rattan in \cite{GR} and by P. \'Sniady in \cite{\'Sn2}. The proof in this article is a new one, which is a consequence of our general combinatorial interpretation.\\

\subsubsection{Low degree terms}

\begin{Th}\label{thlinquad}
The coefficient of the linear monomial $R_{d}$ in $K_k$ is the number of cycles $\tau \in S(k)$ such that $\tau^{-1} \ck$ has $d-1$ cycles. \\

Let $k,j,l$ be positive integers, the coefficient of $R_j R_l$ in $K_k$ is the number (respectively half the number is $j=l$) of pairs $(\tau,\varphi)$ which fulfill the following conditions:
\begin{itemize}
\item The first element $\tau$ is a permutation in $S(k)$ such that $|C(\tau)|=2$. The second element $\varphi$ is a bijection $|C(\tau)| \stackrel{\sim}{\rightarrow} \{1;2\}$. So we count some permutations with numbered cycles.
\item $\tau^{-1} \sigma$ has $j+l-2$ cycles.
\item Among these cycles, at least $j$ have an element in commun with $\varphi^{-1}(1)$ and at least $l$ with $\varphi^{-1}(2)$.
\end{itemize}
\end{Th}
The first part of this theorem was proved by R. Stanley and P. Biane
\cite{Bi2} sepa\-rately, the second is a new result. As in our general
combinatorial interpretation, these
coefficients can be computed by counting permutations in $S(k)$. So, when the support
of the permutations is quite small, we can compute quickly character
values from free cumulants.

\subsection{A combinatorial formula for character values}\label{combap}
The main tool in this article is the following formula\footnote{The
  notations in this article are slightly different with the ones in
  the original papers}, conjectured by R. Stanley in \cite{St2} and
proved by the author in \cite{F\'e}. As noticed in paragraph
\ref{parbackgr}, if we have two sequences $\s{p}$ and $\s{q}$ of
non-negative integers with only finitely many non-zeros terms, we
consider the partition drawn on figure \ref{figpart}:
$$\lambda(\s{p},\s{q}):=\underbrace{\sum_{i \geq 1} q_i, \ldots, \sum_{i \geq 1} q_i}_{p_1 \text{ times}},\underbrace{\sum_{i \geq 2} q_i,\ldots,\sum_{i \geq 2} q_i}_{p_2 \text{ times}},\ldots$$
With this notation, the $R_i(\lambda(\s{p},\s{q}))$ are homogeneous
polynomials of degree $i$ in $\s{p}$ and $\s{q}$.

\begin{Th}\label{thstan}
Let $\s{p}$ and $\s{q}$ be two finite sequences, $\lambda(\s{p},\s{q}) \vdash n$ the associated Young diagram and $\mu \vdash k (k \leq n)$. If $\sigma \in S(k)$ is a permutation of type $\mu$, the character value is given by the formula:
\begin{equation}\label{staneq}
\Sigma_\mu\big(\lambda(\s{p},\s{q})\big)= \sum_{\substack{\dec{\tau} \in S(k)\\ \tau \overline{\tau} = \sigma}} (-1)^{|C(\tau)|+r} N^{\dec{\tau}}(\s{p},\s{q}),
\end{equation}
where $N^{\dec{\tau}}$ is an homogeneous power series of degree $|C(\tau)|$ in $\s{p}$ and $|C(\overline{\tau})|$ in $\s{q}$ which will be defined in section \ref{sect2}.
\end{Th}

This theorem gives a combinatorial interpretation of the coefficients of $\Sigma_\mu$, expressed as a polynomial in variables $\s{p}$ and $\s{q}$. It is natural to wonder if there exists such an expression for free cumulants. Since $R_{l+1}$ is the term of graded degree $l+1$ of $\Sigma_l$ (see \cite[Theorem 1.3]{Bi2}), we have\footnote{A. Rattan has also given a direct proof of this result in \cite{Ra}.}:
\begin{eqnarray}
R_{l+1}(\lambda(\s{p},\s{q})) &=& \sum_{\substack{\dec{\tau} \in S(l)\\ \tau \overline{\tau} = (1 \ldots l) \\ |C(\tau)| + |C(\overline{\tau})| = l+1}} (-1)^{|C(\tau)|+1} N^{\dec{\tau}}(\s{p},\s{q}) ;\nonumber\\
\label{stancum} &=& \sum_{\pi \in NC(l)} (-1)^{|\pi|+1} N^{\pi}(\s{p},\s{q}).
\end{eqnarray}
The second equality comes from the fact that factorizations $\dec{\tau}$ of the long cycle $(1 \ldots l)$ such that $|C(\tau)| + |C(\overline{\tau})| = l+1$ are canonically in bijection with non-crossing partitions (see paragraph \ref{subsectncp}). Note that $N^\pi$ is simply a short notation for $N^{\sigma_\pi,\sigma_\pi^{-1}(1 \ldots l)}$.

From now on, we consider $\Sigma_k$ and $R_l$ as power series in two infinite sets of variables $(\s{p},\s{q})$ and look at equality (\ref{Kerov}) in this algebra (equality as power series in $\s{p}$ and $\s{q}$ is equivalent to equality for all Young diagram $\lambda$, whose size is bigger than a given number). If we expand $K_k(R_2,\ldots,R_{k+1})$, we obtain an algebraic sum of product of power series associated to minimal factorizations. In this article, we write each term of the right side of (\ref{staneq}) as such a sum.

\subsection{Generalized Kerov's polynomials}
The theorems of paragraph \ref{subsectkerpol} corres\-pond to the case where $\mu$ has only one part. But, in fact, they have generalizations for any $\mu \vdash k$.\\

Firstly, there exist universal polynomials $K_\mu$, called generalized Kerov's polynomials, such that:
\begin{equation}\label{kergen}
\Sigma_\mu(\lambda)=K_\mu(R_2(\lambda),\ldots,R_{k+1}(\lambda)).
\end{equation}\\

\begin{eqnarray*}
   \text{Examples: } \Sigma_{2,2} & = & R_3^2 - 4 R_4 - 2 R_2^2 - 2 R_2 ;\\
   \Sigma_{3,2} & = & R_3 \cdot R_4 - 5 R_2 \cdot R_3 - 6 R_5 - 18 R_3 ;\\
   \Sigma_{2,2,2} & = & R_3^3 - 12 R_3 \cdot R_4 - 6 R_3 \cdot R_2^2 + 58 R_3 \cdot R_2 +
   40 R_5 + 80 R_3.
\end{eqnarray*}

Secondly, although these polynomials do not have non-negative coefficients, the following generalization of theorem \ref{mainth} holds:
\begin{Th}\label{thkergen}
Let $\mu \vdash k$ and $\sigma \in S(k)$ a permutation of type $\mu$.
\begin{equation}\label{defkergen}
\Sigma'_{\mu}:= \sum_{\substack{\dec{\tau} \in S(k)\\ \tau \overline{\tau} = \sigma \\ \mathbf{<\dec{\tau}> \text{ trans.}}}} (-1)^{|C(\tau)|+\mathbf{1}} N^{\dec{\tau}},
\end{equation}
where $<\dec{\tau}>$ \emph{trans.} means that the subgroup $<\dec{\tau}>$ of $S(k)$ generated by $\tau$ and $\overline{\tau}$ acts transitively on the set $[k]$.
Then there exists a polynomial $K'_\mu$ with non-negative integer coefficients such that, as power series:
\begin{equation}\label{existkergen}
\Sigma'_\mu=K'_\mu(R_2,\ldots,R_{k+1}).
\end{equation}
\end{Th}

   \begin{eqnarray*}
   \text{Examples: }
   \Sigma'_{2,2} & = & 4 R_4 + 2 R_2^2 + 2 R_2 ;\\
   \Sigma'_{3,2} & = & 6 R_2 \cdot R_3 + 6 R_5 + 18 R_3 ;\\
   \Sigma'_{2,2,2} & = & 64 R_3 \cdot R_2 + 40 R_5 + 80 R_3.
   \end{eqnarray*}

Sections \ref{sect2}, \ref{sect3} and \ref{sect4} are devoted to the proof of this theorem.\\

The quantities $\Sigma'$ are not only practical for the statement of
this theorem, they also appear as disjoint cumulants \cite[Proposition
  22]{F\'S} for study of the asymptotics of character values in
\cite{\'Sn1}. It is also easy to recover $\Sigma$ from $\Sigma'$ by
looking, for each decomposition, at the set partition of $[k]$ in
orbits under the action of $<\dec{\tau}>$ (one has to be careful about the signs):

\begin{equation}\label{sigfromsig2}
\Sigma_{\mu}=\sum_{\Pi \text{ partition of }[l(\mu)]} \left(  \prod_{\{i_1,\ldots,i_l\} \text{ part of } \Pi} (-1)^{l-1} \Sigma'_{\mu_{i_1},\ldots,\mu_{i_l}} \right).
\end{equation}\\

If we invert this formula with (usual) cumulants, then our positivity result on generalized Kerov's polynomials is exactly the one conjectured by A. Rattan and P. \'Sniady in \cite{R\'S}.\\

\subsubsection{Subdominant term for general $\mu$.}
We can also compute some particular coefficients in this general
 context:\\

For low degree terms, the first part of theorem \ref{thlinquad} is still true (it has been proved in \cite{R\'S} in this general context) and the second is true with $K'_\mu$ instead of $K_\mu$  and with an additional condition in the
second part : $<\tau,\tau^{-1} \sigma>$ acts transitively on $[k]$.\\

The highest graded degree in $K'_\mu$ is
$|\mu|+2-l(\mu)$. In the case $l(\mu)=2$, we can explicitly compute
the corresponding term.

\begin{Th}\label{thgenre1gen}
Let $N(l_1,\ldots,l_t;L)$ be the number of solutions of the equation $x_1+\ldots+x_t=L$, fulfilling the condition that, for each $i$, $x_i$ is an integer between $0$ and $l_i$. Then, the coefficient of a monomial $\prod\limits_{i=1}^t R_{j_i}$ of graded degree $r+s$ in $K'_{r,s}$ is:
\begin{equation}\label{eqgenre1gen}
\frac{r \cdot s}{t} \ |\Perm(\s{j})| \ N(j_1-2,\ldots,j_t-2;r-t).
\end{equation}
\end{Th}

This result gives the subdominant term for character values on any fixed permutation:

\begin{Corol}
For any $\mu=(k_1,\ldots,k_r) \vdash k$, one has:
\begin{multline*}\hspace{-1.3 cm} \Sigma_\mu = \prod_{i=1}^r R_{k_i+1} + \sum_{i=1}^r \left[(\prod_{h \neq i} {R_h}) \left(\sum_{|\s{j}|=i-1}  \frac{(k-1)k(k+1)}{24} |\Perm(\s{j})| \prod\limits_i^{l(\s{j})} (j_i-1) R_{j_i} \right) \right] \\
\hspace{-.4cm} + \sum_{1 \leq i_1 < i_2 \leq r} \left[(\prod_{h \neq i_1,i_2} {R_h}) \left(\sum_{|\s{j}|=i_1+i_2} \frac{i_1 \cdot i_2}{l(\s{j})} \ |\Perm(\s{j})| \ N(j_1-2,\ldots,j_t-2;i_1-t) \prod_{i=1}^{l(\s{j})} R_{j_i} \right) \right]\\
+ \text{lower graded degree terms}.
\end{multline*}
\end{Corol}

\begin{proof} In equation (\ref{sigfromsig2}), the only summands which contain terms of degree $|\mu|+r-2$ are the one indexed by the partition of $[l(\mu)]$ in singletons and those indexed by partitions in one pair and singletons.\end{proof}

\subsection{Organization of the article}
In section \ref{sect2}, we will associate a map to each pair of
permutations. This will help us to define the associated power series
$N$. In section \ref{sect3}, for any map $M$, we write $N(M)$ as an
algebraic sum of power series associated to minimal
factorizations. The section \ref{sect4} is the end of the proof of
theorem \ref{thkergen}. Then, in section \ref{sect5}, we will compute
some particular coefficients (proofs of theorems \ref{thgenre1}, \ref{thlinquad} and \ref{thgenre1gen}).

\section{Maps and polynomials}\label{sect2}
In this section, we define the power series $N^{\dec{\tau}}$ as the composition of three functions:
$$\xymatrix{S(k) \times S(k) \ar@{->}[r]^{\S \ \ref{permtomaps} \qquad} & \text{bicolored labeled map} \ar@{->}[r]^{\quad \text{Forget}} & \text{bicolored graph} \ar@{->}[r]^{\quad \S \  \ref{graphtops}} & \C[[\textbf{p},\textbf{q}]]}$$

\subsection{From permutations to maps}\label{permtomaps}
Let us give some definitions about graphs and maps.
\begin{Def}[graphs]\label{defgr}
\begin{itemize}
\item A graph is given by:
\begin{itemize}
\item a finite set of vertices $V$ ;
\item a set of half-edges $H$ with a map \emph{ext} from $H$ to $V$
  (the image of an half-edge is called its extremity) ;
\item a partition of $H$ into pairs (called edges, whose set is denoted $E$) and
  singletons (the external half-edges).
\end{itemize}
\item A bicolored graph is a graph with a partition of $V$ in two sets (the set of white vertices $V_w$ and the set of black vertices $V_b$) such that, for each edge, among the extremities of its two half-edges, one is black and one is white.
\item A labeled graph is a graph with a map $\iota$ from $E$ in $\N^\star$. Moreover, we say that it is well labeled if $\iota$ is a bijection of image $[|E|]$.
\item An oriented edge $ {e}$ is an edge $e$ with an order of its two half-edges.
\item An oriented loop is a sequence of oriented edge $ {e_1},\ldots, {e_l}$ such that:
\begin{itemize}
\item For each $i$, the extremity $v_i$ of the first half-edge of $ {e_{i+1}}$ is the same as the extremity of the second of $ {e_i}$ (with the convention $e_{l+1}=e_1$);
\item All the $v_i$'s and the $e_i$'s are different (an edge does not appear twice, even with different orientations).
\end{itemize}
We identify sequences that differ only by a cyclic permutation of their orien\-ted edges.
\item The free abelian group on graphs has a natural ring structure: the product of two graphs is by definition their disjoint union.
\end{itemize}
\end{Def}
\begin{Def}[Maps]\label{defmaps}
\begin{itemize}
\item A map is a graph supplied with, for each vertex $v$, a cyclic order on the set of all half-edges (including the external ones) of extremity $v$ (\textit{i.e.} $\text{ext}^{-1}(v)$).
\item Consider an half-edge $h$ of a map $M$. Thanks to the map structure, there is a cyclic order on the set of half-edges having the same extremity as $h$. We call successor of $h$ the element just after $h$ in this order.
\item Since a map is a graph with additional informations, we have the notion of bicolored and/or (well-)labeled map.
\item A face of a map is a sequence of oriented edge $ {e_1},\ldots, {e_k}$ such that, for each $i$, the first half-edge of $e_{i+1}$ ($e_{l+1}=e_1$) is the successor of the second half-edge of $e_i$. As for loops, we identify the sequences which differ by cyclic permutations of their orien\-ted edges. Then each oriented edge is in exactly one face.
\item If $F$ is a face of a map is labeled and bicolored, we denote by $E(F)$ the set of edges appearing in $F$ with the white to black orientation. The word associated to a face is the word $w(F)$ of the labels of the elements of $E(F)$ (it is defined up to a cyclic permutation).
\item A face, which is also a loop (all vertices and edges of the face are distinct) and which does not contain an external half-edge, is called a polygon.
\end{itemize}
\end{Def}

\begin{Rem}
A map, whose underlying graph is a tree, is a planar tree. It has exactly one face.
\end{Rem}

\begin{figure}
\begin{center}
\includegraphics{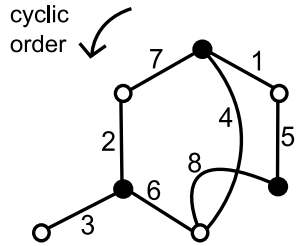}\caption{\footnotesize Example of a bicolored labeled map, with exactly one face whose associated word is $12345678$}
\label{figexcarte}
\end{center}
\end{figure}

\subsubsection{Map associated with a pair of permutations}
The following construction is classical (it generalizes the work of I.P. Goulden and D.M. Jackson in \cite{GJ}) but we recall it for completeness.
\begin{Def}
To a well-labeled bicolored map $M$ with $k$ edges and no external half-edges, we associate the pair of permutations $(\dec{\tau}) \in S(k)^2$ defined by: if $i$ is an integer in $[k]$, $e$ the edge of $M$ with label $i$ and $h$ its half-edge with a white (resp. black) extremity, then $\tau(i)$ (resp. $\overline{\tau}(i)$) is the label of the edge containing the successor of $h$.
\end{Def}

It is easy to see that this defines a bijection between well-labeled bicolored maps and pairs of permutations in $S(k)$. Its inverse associates to a pair of permutations $(\dec{\tau})$ the following bicolored labeled map $M^{\dec{\tau}}$: the set of white vertices is $C(\tau)$, the one of black vertices $C(\overline{\tau})$, the set of half-edges $\{1^w,1^b,\ldots,k^w,k^b\}$ is partitioned in edges $\{i^w,i^b\}$ and the cycle $(i_1,\ldots,i_l)$ of $\tau$ (resp. $(j_1,\ldots,j_l)$ of $\overline{\tau}$) is the extremity of the half-edges $i_1^w,\ldots,i_l^w$ (resp. $j_1^b,\ldots,j_l^b$) in this cyclic order.

The following property follows straight forward from the definition:
\begin{Prop}\label{prodface}
The words associated to the faces of $M^{\dec{\tau}}$ are exactly the cycles of the product $\tau \overline{\tau}$.
\end{Prop}

\begin{example}The map drawn on figure \ref{figexcarte} is associated to the pair of permutations $\big((15)(27)(3)(486),(174)(236)(58)\big)$ of product $(12345678)$. The word associated to its unique face is $12345678$ as predicted by proposition \ref{prodface}.\\
\end{example}

Note that the connected components of $M^{\dec{\tau}}$ are in bijection with the orbits of $[k]$ under the action of $<\dec{\tau}>$. So, a factorization is transitive if and only if its map is connected. In particular, maps of minimal factorizations of the full cycle $\ck$ are exactly the connected maps with $k+1$ vertices and $k$ edges, that is to say the planar trees.

\subsection{From graphs to polynomials}\label{graphtops}
\begin{Def}
Let $G$ be a bicolored graph and $V$ its set of vertices, disjoint union of $V_b$ and $V_w$. An evaluation $\psi: V \rightarrow \N^\star$ is said admissible if, for any edge between a white vertex $w$ and a black one $b$, it fulfills $\psi(b) \geq \psi(w)$. The power series $N(G)$ in indeterminates $\s{p}$ and $\s{q}$ is defined by the formula:
\begin{equation}\label{defF}
N(G) = \sum_{\substack{\psi: V \rightarrow \N \\ \text{admissible} }} \prod_{w \in V_w} p_{\psi(w)} \prod_{b \in V_b} q_{\psi(b)}.
\end{equation}
Note that $N$ is extended to the ring $\mathbb{A}_{bg}$ of bicolored graphs by $\Z$-linearity. It is in fact a morphism of rings (the power series associated to a disjoint union of graphs is simply the product of the power series associated to these graphs).\\
If $\tau$ and $\overline{\tau}$ are two permutations in $S(k)$, we put:
$$N^{\dec{\tau}}:=N(M^{\dec{\tau}}).$$
\end{Def}

This definition is the one that appears in theorem \ref{thstan}. The main step of our proof of Kerov's conjecture is to write the power series associated to any pair of permutations as an algebraic sum of power series associated to forests (\textit{i.e.} pro\-ducts of power series associated to minimal factorizations).\\

\begin{figure}
\begin{center}
\includegraphics{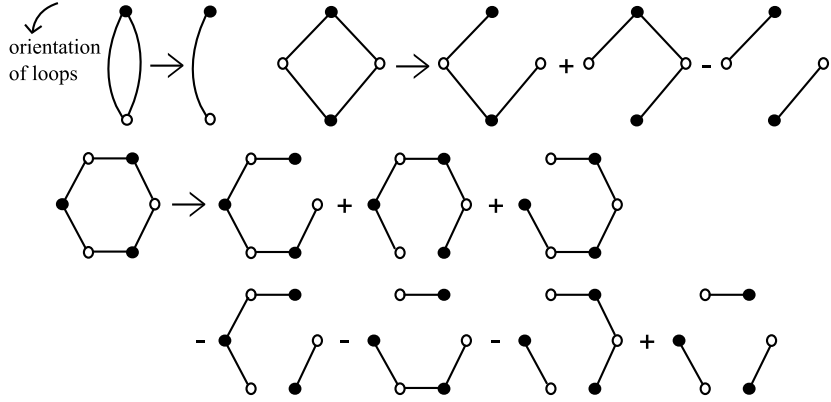}\caption{\footnotesize Illustration of definition of transformation $T_{{L}}$}
\label{figtransel}
\end{center}
\end{figure}

Let $G$ be a bicolored graph and $L$ an oriented loop of $G$. We denote by $E(L)$ the set of edges which appear in the sequence $L$ oriented from their white extremity to their black one. Let us define the following element of the $\Z$-module $\mathbb{A}_{bg}$:
\begin{equation}\label{transel}
T_{{L}}(G)=\sum_{\substack{E' \subset E(L)\\E' \neq \emptyset}} (-1)^{|E'|-1} G \backslash E',
\end{equation}
where $G \backslash E'$ denotes the graph obtained by taking $G$ and erasing its edges belonging to $E'$ (it is a subgraph of $G$ with the same set of vertices). These elementary transformations are drawn on figure \ref{figtransel}, where we have only drawn vertices and edges belonging to the loop ${L}$ (so these schemes can be understood as local transformations).
\\

\begin{figure}
\begin{center}
\includegraphics{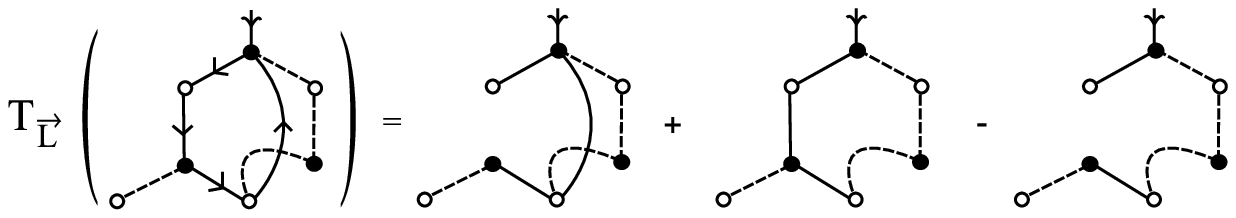}\caption{\footnotesize Example of an elementary transformation.}
\label{figextrel}
\end{center}
\end{figure}

An example of such a transformation is drawn in figure \ref{figextrel}. $G$ is the map of figure $\ref{figexcarte}$ (we forget the labels and the map structure) and ${L}$ the loop $7,2,6,4$.\\

We have the following conservation property:
\begin{Prop}\label{propinv}
If $G$ is a bicolored graph and ${L}$ an oriented loop of $G$, then
\begin{equation}\label{transFinv}
N\big(T_{{L}}(G)\big)=N(G).
\end{equation}
\end{Prop}

\begin{proof}
Let $G$ be a bicolored graph and $V_w$, $V_b$, $E$ as in definition \ref{defgr}. We write the series $N(G)$ as the following sum:
\begin{eqnarray}
N(G) & = & \sum_{\psi_w: V_w \rightarrow \N^\star} \left[\sum_{\substack{\psi: V \rightarrow \N^\star \text{admissible} \\ \psi_{/V_w}=\psi_w}} \prod_{w \in V_w} p_{\psi(w)} \prod_{b \in V_b} q_{\psi(b)} \right] ;\nonumber\\
\label{fassum} & = & \sum_{\psi_w: V_w \rightarrow \N^\star} N_{\psi_w}(G).
\end{eqnarray}

Since all the graphs in the equality (\ref{transFinv}) have the same set of vertices $V_w$, it is enough to prove that, for every $\psi_w: V_w \rightarrow \N^\star$, we have:
\begin{equation}\label{fpsi}
N_{\psi_w} \big(T_{{L}}(G)\big)=N_{\psi_w}(G).
\end{equation}

Let us fix a partial evaluation $\psi_w: V_w \rightarrow \N^\star$. If we choose a numbering $w_1,\ldots,w_l$ (with respect to the loop order) of the white vertices of ${L}$, then there exists an index $i$ such that $\psi_w(w_{i+1}) \geq \psi_w(w_i)$ (with the convention $w_{l+1}=w_1$). Denote by $e$ the edge just after $w_i$ in the loop ${L}$. It is an erasable edge. So we have a bijection:
\begin{eqnarray*}
\big\{E' \subset E(L), e \notin E'\big\} & \stackrel{\sim}{\rightarrow} & \big\{E'' \subset E(L), e \in E''\big\}\\
E' & \mapsto & E''=E' \cup \{e\}.
\end{eqnarray*}

But, this bijection has the following property:
$$N_{\psi_w}(G \backslash E')=N_{\psi_w}(G \backslash (E' \cup \{e\}).$$
Indeed the admissible evaluations whose restrictions to white vertices is $\psi_w$ are the same for the two graphs $G \backslash E'$ and $G \backslash (E' \cup \{e\})$. The only thing to prove is that, if such a $\psi$ is admissible for $G \backslash (E' \cup \{e\})$, it fulfills also: $\psi(b_e) \geq \psi(w_i)$, where $b_e$ is the black extremity of $e$. This is true because $$\psi(b_e) \geq \psi(w_{i+1}) = \psi_w(w_{i+1}) \geq \psi_w(w_i) = \psi(w_i).$$

To conclude the proof, note that cardinals of $E'$ and $E' \cup \{e\}$ have different parity so they appear with different signs in $G - T_{{L}}(G)$. Their contributions to (\ref{fpsi}) cancel each other and the proof is over.
\end{proof}

Recall that $N$ is a morphism of rings, so $(\mathbb{A}_{bg})_{/Ker N}$ is a ring.
\begin{Corol}\label{arbeng}
The ring $(\mathbb{A}_{bg})_{/Ker N}$ is generated by trees.
\end{Corol}
\begin{proof} Just iterate the proposition by choosing any oriented loop until there is no loop left (if a graph is not a disjoint union of trees, there is always one).
\end{proof}

However, forests are not linearly independent in $(\mathbb{A}_{bg})_{/Ker N}$.

\section{Map decomposition}\label{sect3}
By iterating proposition \ref{propinv} until there are only forests
left, given a graph $G$, we obtain an algebraic sum of forests whose
associated power series is $N(G)$. But there are many possible choices of oriented loops and they can give different sums of forests. In this section, we explain, how, by restricting the choices, we choose a particular one, which depends on the map structure and the labeling.

\subsection{Elementary decomposition}
To do coherent choices, it is convenient to add an external half-edge to our map. So, in this paragraph, we deal with bicolored maps with exactly one external half-edge $h$. They generate a free $\Z$-module denoted $\mathbb{A}_{bm,1}$.\\

If $M$ is such a map, let $\star$ be the extremity of its external half-edge. An (oriented) loop ${L}$ is said admissible if:
\begin{itemize}
\item The vertex $\star$ is a vertex of the loop, that is to say that $\star$ is the extremity of the second half-edge $h_{i,2}$ of $e_i$ and of the first half-edge $h_{i+1,1}$ of $e_{i+1}$ for some $i$;
\item The cyclic order at $\star$ restricted to the set $\big\{h,h_{i,2},h_{i+1,1}\big\}$ is the cyclic order $\big( h,h_{i+1,1},h_{i,2} \big)$.
\end{itemize}
For example, the oriented loop ${L}$ of figure \ref{figextrel} is admissible. If ${L}$ satisfies the first condition, exactly one among the oriented loops ${L}$ and ${L}'$ is admissible (where ${L}'$ is ${L}$ with the opposite orientation).

\begin{Defth}\label{invdecel}
There exists a unique linear operator $$D_1: \mathbb{A}_{bm,1} \rightarrow \mathbb{A}_{bm,1}$$ such that:
\begin{itemize}
\item The image of a given map $M$ lives in the vector space spanned by its submaps with the same set of vertices ;
\item If ${L}$ is an admissible loop of $M$, then
\begin{equation}\label{defD1}
D_1(M)=D_1 \big(T_{{L}}(M) \big).
\end{equation}
Note that this equality is meant as an equality between submaps of $M$, not only as isomorphic maps ;
\item If there is no admissible loops in $M$, then $D_1(M)=M$.
\end{itemize}
\end{Defth}

\begin{proof}
If $M$ is a bicolored map, all graphs appearing in $T_{{L}}(M)$ have stricly less edges than $M$. So the uniqueness of $D_1$ is obvious.\\

The existence of $D_1$ will be proved by induction. Denote, for every $N$, $\mathbb{A}^N_{bm,1}$ the submodule of $\mathbb{A}_{bm,1}$ generated by graphs with at most $N$ edges. We will prove that there exists, for every $N$, an operator $D^N_1: \mathbb{A}^N_{bm,1} \rightarrow \mathbb{A}^N_{bm,1}$, extending $D^{N-1}_1$ if $N \geq 1$, and satisfying the conditions asked for $D_1$. The case $N=0$ is very easy because $\mathbb{A}^0_{bm,1}$ is generated by graphs without admissible loops, so $D^0_1=Id$. If our statement is proved for any $N$, it implies the existence of $D_1$: take the inductive limit of the $D^N_1$.\\

Let $N \geq 1$ and suppose that $D^{N-1}_1$ has been built. To prove the existence of $D^N_1$, we have to prove that, if $M$ has admissible loops, then $D^{N-1}_1 \big(T_{{L}}(M) \big)$ does not depend on the chosen admissible loop ${L}$. To do this, let us denote by $M_\star$ the submap of $M$ containing exactly all the edges of $M$ which belong to some admissible loop of $M$. The maps $M$ and $M_\star$ have exactly the same admissible loops. We define $H=|E(M_\star)|-|V(M_\star)|+1$ (which might be understood as the number of independent loops in $M_\star$).\\

If $H=0,1$, the map $M$ has at most one admissible loop, so there is nothing to prove:
\begin{itemize}
\item If $M$ has exactly no admissible loop, then $D^N_1(M)=M$.
\item If $M$ has exactly one admissible loop ${L}$, then $D^N_1(M)=T_{{L}}(M)$.
\end{itemize}

If $H=2$ and if there is a vertex of valence $4$ in $M_\star$ different from $\star$, then there is at most one admissible loop. If H=2 and if $\star$ is a vertex of valence $4$, then there are two admissible loops ${L_1}$ and ${L_2}$ without any edges in commun, so the transformation with respect to these loops commute, so
$$D^{N-1}\big(T_{{L_1}}(M)\big) = T_{L_2}\big(T_{{L_1}}(M)\big) =T_{L_1}\big(T_{{L_2}}(M)\big) = D^{N-1}\big(T_{{L_2}}(M)\big).$$

If $H=2$ and if $\star$ and an other vertex $v$ have valence $3$, there are three admissible loops. In $M_\star$, there are three different paths $c_0,c_1,c_2$ going (without any repetition of vertices or edges) from $\star$ to $v$. We number them such that, if $h_i$ is the first half-edge of the path $c_i$, the cyclic order at $\star$ is $(h,h_0,h_1,h_2)$. Let us denote by $E_i (0 \leq i \leq 2)$ (resp. by $E_{\bar{i}}$) the set of edges appearing in $c_i$ oriented from their black vertex to their white one (resp. from their white vertex to their black one). If $I=\{i_1,\ldots,i_l\} \subset \{0,1,2,\bar{0},{\bar{1}},{\bar{2}}\}$, we consider the following element of $\mathbb{A}_{bg,1}$:
$$M_I = \sum_{\emptyset \neq E'_1 \subset E_{i_1},\ldots,  \emptyset \neq E'_l \subset E_{i_l}} (-1)^{|E'_1|-1} \ldots (-1)^{|E'_l|-1} M \backslash \big( E'_1 \cup \ldots \cup E'_l \big).$$
Let ${L_1}=c_0 \cdot \overline{c_1}$, ${L_2}=c_1 \cdot \overline{c_2}$ and ${L_3}=c_0 \cdot \overline{c_2}$ be the three admissible loops of $M$. Their respective sets of erasable edges are $E_{\bar{0}} \cup E_1$, $E_{\bar{0}} \cup E_2$ and $E_{\bar{1}} \cup E_2$. So we have (the figure \ref{figexinv} shows this computation on an example, where all sets $E_i$ are of cardinal $1$):
\begin{figure}
\begin{center}
\includegraphics{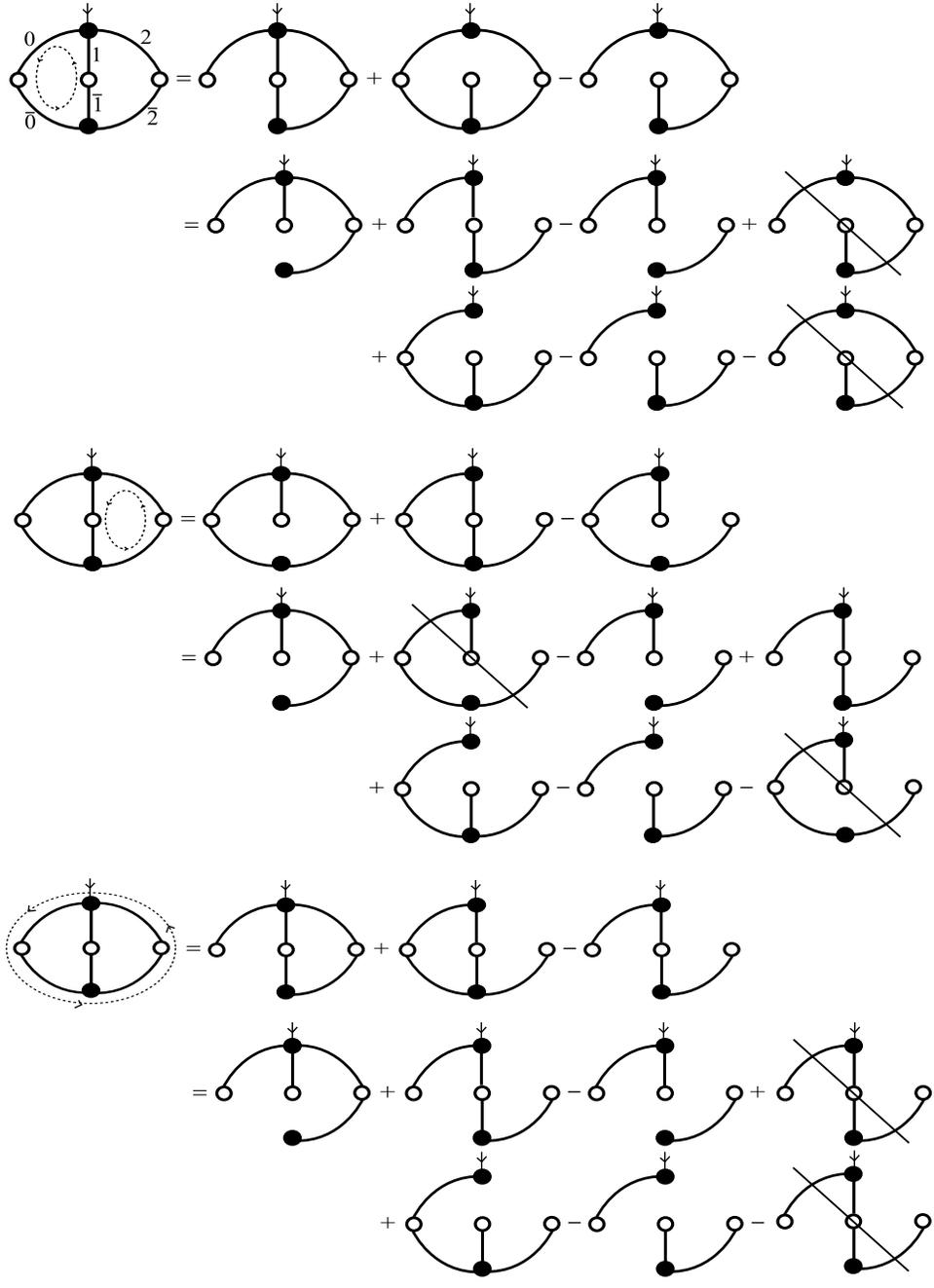}\caption{\footnotesize One particular case of definition-theorem \ref{invdecel}}
\label{figexinv}
\end{center}
\end{figure}
\begin{eqnarray*}
T_{{L_1}}(M) &=& \sum_{\substack{E' \subset E_1 \\ E' \neq \emptyset}} (-1)^{|E'|-1} M \backslash E' + \sum_{\substack{E' \subset E_{\bar{0}} \\ E' \neq \emptyset}} (-1)^{|E'|-1} M \backslash E' \\ & & + \sum_{\substack{E' \subset (E_1 \cup E_{\bar{0}})\\(E' \cap E_1) \neq \emptyset, (E' \cap E_{\bar{0}}) \neq \emptyset}} (-1)^{|E'|-1} M \backslash E' ;\\
&=& M_{\bar{0}} + M_1 - M_{1,{\bar{0}}}.
\end{eqnarray*}
For each graph appearing in $M_{\bar{0}}$, $M_1$  there is only one admissible loop so $D^{N-1}_1$ is just given by the corresponding elementary transform:
\begin{eqnarray*}
D^{N-1}_1(T_{{L_1}}(M)) &=& M_{{\bar{0}},{\bar{1}}} + M_{2,{\bar{0}}} - M_{2,{\bar{0}},{\bar{1}}} + M_{1,{\bar{0}}} + M_{1,2} - M_{1,2,{\bar{0}}} -M_{1,{\bar{0}}},\\
&=& M_{{\bar{0}},{\bar{1}}} + M_{2,{\bar{0}}} - M_{2,{\bar{0}},{\bar{1}}} + M_{1,2} - M_{1,2,{\bar{0}}}.
\end{eqnarray*}
For the other admissible loops, we obtain:
\begin{eqnarray*}
D^{N-1}_1(T_{{L_2}}(M)) &=& D^{N-1}_1(M_{\bar{1}} + M_2 - M_{2,{\bar{1}}}),\\
&=&  M_{{\bar{0}},{\bar{1}}} + M_{2,{\bar{1}}} - M_{2,{\bar{0}},{\bar{1}}} + M_{2,{\bar{0}}} + M_{1,2} - M_{1,2,{\bar{0}}} -M_{2,{\bar{1}}},\\
&=&  M_{{\bar{0}},{\bar{1}}} - M_{2,{\bar{0}},{\bar{1}}} + M_{2,{\bar{0}}} + M_{1,2} - M_{1,2,{\bar{0}}} ;\\
D^{N-1}_1(T_{{L_3}}(M)) &=& D^{N-1}_1(M_{\bar{0}} + M_2 - M_{2,{\bar{0}}}),\\
&=&  M_{{\bar{0}},{\bar{1}}} + M_{2,{\bar{0}}} - M_{2,{\bar{0}},{\bar{1}}} + M_{2,{\bar{0}}} + M_{1,2} - M_{1,2,{\bar{0}}} -M_{2,{\bar{0}}},\\
&=&  M_{{\bar{0}},{\bar{1}}} - M_{2,{\bar{0}},{\bar{1}}} + M_{2,{\bar{0}}} + M_{1,2} - M_{1,2,{\bar{0}}}.\\
\end{eqnarray*}

If $H=2$ and if there are two vertices $v$ and $v'$ of valence 3 distinct from $\star$, the proof is similar. We use the same notations, except that:
\begin{itemize}
\item The paths $c_0$, $c_1$ and $c_2$ go from $v$ to $v'$.
\item The vertex $\star$ is on $c_0$. It does not matter to exchange $c_1$ and $c_2$.
\item If the half-edge just before (resp. just after) $\star$ in $c_0$ is denoted by $h_1$ (resp. $h_2$), the cyclic order at $\star$ induces the order $(h_1,h,h_2)$.
\end{itemize}
In this case, there are only two admissible loops ${L_1}$ and ${L_3}$ in $M$ and a little computation proves the theorem:
\begin{eqnarray*}
D^{N-1}_1(T_{{L_1}}(M)) &=& D^{N-1}_1(M_{\bar{0}} + M_1 - M_{1,{\bar{0}}}),\\
&=& M_{\bar{0}} + M_{1,{\bar{0}}} + M_{1,2} - M_{1,2,{\bar{0}}} -M_{1,{\bar{0}}},\\
&=& M_{\bar{0}} + M_{1,2} - M_{1,2,{\bar{0}}};\\
D^{N-1}_1(T_{L_{\bar{0}}}(M)) &=& D^{N-1}_1(M_{\bar{0}} + M_2 - M_{2,{\bar{0}}}),\\
&=&  M_{\bar{0}} + M_{2,{\bar{0}}} + M_{1,2} - M_{1,2,{\bar{0}}} -M_{2,{\bar{0}}},\\
&=&  M_{\bar{0}} + M_{1,2} - M_{1,2,{\bar{0}}}.
\end{eqnarray*}
The proof is over in the case $H=2$.\\

The case $H \geq 3$ needs the two following lemmas:

\begin{Lem}\label{lemsanse}
Let ${L}$ be an admissible loop of $M$ and $e$ an edge of $M \backslash {L}$. Then,
$$D^{N-1}_1\big(T_{{L}}(M)\big) = D^{N-1}_1\big(D^{N-1}_1(M\backslash \{e\}) \cup \{e\} \big),$$
where, for a submap $M' \subset M$ with the same set of vertices which does not contain $e$, $M' \cup \{e\}$ the map obtained by adding the edge $e$ to $M'$.
\end{Lem}

\begin{proof}
To compute the left side of the equation, we choose, for every graph in $T_{{L}}(M)$, one of its admissible loop, apply the associated transformation and iterate this. If, whenever it is possible, we choose an admissible loop that does not contain $e$, the first choices done are also choices of admissible loops for the map $M\backslash \{e\}$. After the associated transformations, we obtain $D^{N-1}_1(M\backslash \{e\}) \cup \{e\}$ and the lemma follows.
\end{proof}

\begin{Lem}\label{lemtechboucl}
If $H \geq 3$ and if ${L_1}$ and ${L_2}$ are two admissible loops with ${L_1} \cup {L_2} = M$, then there exists a third one ${L}$ such that ${L} \cup {L_1} \neq M$ and ${L} \cup {L_2} \neq M$.
\end{Lem}

\begin{proof}
We choose a numbering of the oriented edges of the loops so that the first half-edge of $e_1$ has $\star$ for extremity. We suppose (eventually by exchanging ${L_1}$ and ${L_2}$) that the first half-edge of ${L_1}$ is between $h_0$ and the first half-edge of ${L_2}$ in the cyclic order of $\star$. As ${L_1} \cup {L_2} = M$, the loops ${L_1}$ and ${L_2}$ have an other vertex in commun than $\star$ (otherwise, $M$ is a wedge of two cycles and $H=2$). Let $v$ be the first vertex of ${L_1}$ which is also in ${L_2}$ but such that the paths from $\star$ to $v$ given by the beginnings of ${L_1}$ and ${L_2}$ are different. Let us consider the sequence ${L}$ equal to the concatenation of the beginning of ${L_1}$ (from $\star$ to $v$) and the end of ${L_2}$ (from $v$ to $\star$). With this definition:
\begin{itemize}
\item All vertices and edges appearing in ${L}$ are distinct. Moreover, ${L}$ is an admissible loop ;
\item The edge before $v$ in ${L_2}$ belongs neither to ${L_1}$ nor to ${L}$ ;
\item As $H >2$, the ends of ${L_1}$ and ${L_2}$ (from $v$ to $\star$) are different. So there is an edge in the end of ${L_1}$ which belongs neither to ${L_2}$ nor to ${L}$.\qedhere
\end{itemize}
\end{proof}

Lemma \ref{lemsanse} implies: if ${L_1}$ and ${L_2}$ are admissible loops such that ${L_1} \cup {L_2} \neq M$, we have:
$$D_1\big(T_{{L_1}}(M)\big) = D_1\big(T_{{L_2}}(M)\big).$$
Together with lemma \ref{lemtechboucl}, this ends the proof of the theorem.
\end{proof}

\begin{Rem}[useful in paragraph \ref{cummapdec}]\label{decnonbic}
The definition of this operator does not really need the maps to be bicolored. It is enough to suppose that each edge has a privileged orientation. In this context, the erasable edges of a oriented loop are the one which appear in the loop in their privileged orientation and operator $T_{{L}}$ has a sense. A bicolored map can be seen this way if we choose as orientation of each edge the one from the white vertex to the black one.
\end{Rem}

\subsection{Complete decomposition}\label{comdec}

\begin{figure}
\begin{center}
\includegraphics{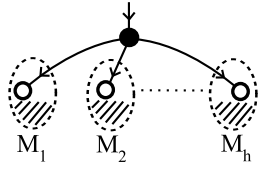}\caption{\footnotesize General form of the connected component containing $\star$ of a map appearing in $D_1(M)$.}
\label{figsortieD1}
\end{center}
\end{figure}

It is immediate from the definition that every map $M'$ appearing with a non-zero coefficient in $D_1(M)$ has no admissible loops. Thus they are of the following form (drawn on figure \ref{figsortieD1}):\\
The vertex $\star$ is the extremity of half-edges $h_i (0 \leq i \leq l)$, including the external one $h_0$, numbered with respect to the cyclic order. For $i \geq 1$, $h_i$ belongs to an edge $e_i$, whose other extremity is $v_i$. Each $v_i$ is in a different connected component $M_i$ (called leg) of $M \backslash \{h_1,\ldots,h_l\}$. Note that we have only erased the half-edge $h_i$ and not the whole edge $e_i$ so that each $M_i$ keeps an external half-edge.\\

If we have a family of submaps $M'_i=M_i \backslash \{E'_i\}$ of the $M_i$ we consider the map $\phi_M(M'_1,\ldots,M'_l)=M \backslash \bigcup \{E'_i\}$ obtained by replacing in $M$ each $M_i$ by $M'_i$.\\

The outcome of operator $D_1$ is an algebraic sum of maps, which are much more complicated than planar forests. So, in order to write $N(M)$ as an algebraic sum of series associated to minimal factorizations, we have to iterate such operations.\\

We want to define decompositions of maps associated to pairs of permutations, so of well-labeled bicolored maps without external edges. But it is convenient to work on a bigger module: the ring $\mathbb{A}_{blm,\leq 1}$ of bicolored labeled maps with at most one external half-edge per connected component.

\begin{Defprop}\label{defD}
There exists a unique linear operator $$D: \mathbb{A}_{blm,\leq 1} \rightarrow \mathbb{A}_{blm,\leq 1}$$ such that:
\begin{enumerate}
\item If $M$ has only one vertex, then $D(M)=M$ ;
\item If $M$ has more than one connected components $M=\prod M_i$, then one has $D(M)=\prod D(M_i)$ ;
\item If $M$ has only one connected component and no external half-edge, consider its edge $e$ of smallest label. Let $h$ be the half-edge of $e$ of black extremity. We denote by $\overline{M}$ the map obtained by adding one external half-edge between $h$ and its successor. Then $D(M)=D(\overline{M})$ ;
\item If $M$ has only one connected component with one half-edge but no admissible loops, we use the notations of the previous paragraph. As the $M_i$ are connected maps with an external half-edge, we can compute $D(M_i)$ (third or fifth case). Then $D(M)$ is given by the formula:
$$D(M)=\phi_M(D(M_1),\ldots,D(M_l)),$$
where $\phi_M$ is extended by multilinearity to algebraic sums of submaps of the $M_i$'s.
\item Else, $D(M)=D(D_1(M))$.
\end{enumerate}
\end{Defprop}

Existence and uniqueness of $D$ are obvious. The image of a map $M$ by $D$ is in the subspace generated by its submaps with the same set of vertices, no isolated vertices and no loops, \text{i.e.} its covering forests without trivial trees. Note also that forests are fixed points for $D$ (immediate induction).\\

\begin{example}
We will compute $D(M)$ where $M$ is the map of the figure \ref{figcarteadec} (without the external half-edge).\\
\begin{figure}
\begin{center}
\includegraphics{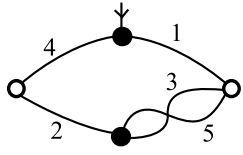}\caption{\footnotesize Map $\overline{M}$.}
\label{figcarteadec}
\end{center}
\end{figure}
The map $M$ belongs to the third kind, so we have to add an external half-edge as on the figure. Now, $\overline{M}$ is a map of the fifth type and we have to compute $D_1(\overline{M})$: this is very easy because the two transformations associated with admissible loops lead to the same sum of submaps which do not contain any admissible loop.
\begin{eqnarray*}
D_1(M) &=& M \backslash \{1\} + M \backslash \{2\} - M \backslash \{1,2\};\\
\text{So }D(M) &=& D(M \backslash \{1\}) + D(M \backslash \{2\}) - D(M \backslash \{1,2\}).
\end{eqnarray*}
The map $M \backslash \{1\}$ is a map of the fourth type with only one leg $M_1$, which is drawn at figure \ref{figm1234}.\\
\begin{figure}
\begin{center}
\includegraphics{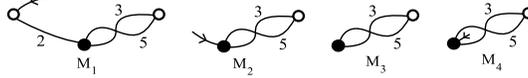}\caption{\footnotesize Maps involved in the computation of the example.}
\label{figm1234}
\end{center}
\end{figure}
This map $M_1$ is again of the fourth type (with one leg: the map $M_2$ of the figure \ref{figm1234}) so we have to compute $D(M_2)$, which is simply $D_1(M_2)=M_2 \backslash \{5\}$. This implies immediately that $D(M_1)=M_1 \backslash \{5\}$ and:
\begin{eqnarray*}
D(M \backslash \{1\}) &=& M \backslash \{1,5\}.\\
\text{Similarly, }D(M \backslash \{2\}) &=& M \backslash \{2,3\}.
\end{eqnarray*}
Now we look at the map $M \backslash \{1,2\}$. It has two connected component (we have to apply rule 2): one is a tree and has a trivial image by $D$, the other one $M_3$ has no external half-edge. We have to add one external half-edge to $M_3$ with the third rule and obtain $M_4$. Now, it is clear that $D_1(M_4) = M_4 \backslash \{3\}$, so one has $D(M \backslash \{1,2\})=M \backslash \{1,2,3\}$.\\
Finally
$$D(M)=M \backslash \{1,5\}+M \backslash \{2,3\}-M \backslash \{1,2,3\}.$$
\end{example}

As we can see on the example, when we replace $M_i$ by its image by several elementary transformations in $M$, we obtain the image of $M$ by the same transformations. So, by an immediate induction, the operator $D$ consists in applying to $M$ an elementary transformation $T_{{L}}$ (with restricted choices), then one to each map of the result which is not a forest, \textit{etc.} until there are only forests left. An immediate consequence is the $D-$invariance of $N$.\\

\begin{Rem}\label{ordre}
Note that transformations indexed by loops which are in different connected components and/or in different legs of the map (fourth case) commute.
\end{Rem}

\subsection{Signs}
In this paragraph, we study the sign of the coefficients in the expression $D(M)$. This is central in the proof of theorem \ref{thkergen} because we will show that the coefficients of $K'_\mu$ can be written as a sum of coefficients of $D(M)$, for some particular maps $M$.

\begin{Prop}\label{propsign}
Let $M' \subset M$ two maps with the same set of vertices and respectively $t_{M'}$ and $t_M$ connected components. The sign of the coefficient of $M'$ in $(-1)^{t_{M}}D(M)$ is $(-1)^{t_{M'}}$.
\end{Prop}
\begin{proof}
Due to the inductive definition of $D$ using $D_1$, it is enough to prove the result for operator $D_1$ in the case where $M$ is a connected ($t_M=1$) bicolored map with one external half-edge. We proceed by induction over the number of edges in $M \backslash M'$. If $M'=M$, the result is obvious. Note that if $M'$ has a non-zero coefficient in $D_1(M)$, we have necessarily $M \backslash M' = \{e_1,\ldots,e_l\}$ where each $e_i$ belongs at least to one admissible loop.\\

\emph{First case:} There exists an edge $e \in M \backslash M'$ such that $M \backslash \{e\}$ has at least one admissible loop. Let us define $M_1=M \backslash \{e\}$ and apply the lemma \ref{lemsanse}: $D_1(M)=D_1\big(D_1(M_1) \cup \{e\}\big)$. The submaps $M''$ of $M_1$ containing $M'$ can be divided in two classes:
\begin{itemize}
\item Either $M'' \cup \{e\}$ has the same number $t$ of connected components as $M''$. By induction hypothesis, the sign of the coefficient of $M'' \cup \{e\}$ in $D_1(M_1) \cup \{e\}$ is $(-1)^{t-1}$ ;
\item Or $M'' \cup \{e\}$ has strictly less connected components than $M''$. In this case $\{e\}$ does not belong to any loops of $M'' \cup \{e\}$, so every graph appearing in $D_1(M'' \cup \{e\})$ does contain $\{e\}$. In particular, the coefficient of $M'$ in $D_1(M'' \cup \{e\})$ is zero.
\end{itemize}
Finally, the coefficient of $M'$ in $D_1(M)$ is the same as in the sum of $D_1(M'' \cup \{e\})$ for $M''$ of the first class. So the result comes from the induction hypothesis applied to $M' \subset M'' \cup \{e\}$ (which can be done because $M'' \cup \{e\}$ has strictly less edges than $M$).\\

\emph{Second case:} Else, up to a new numbering of edges of $M \backslash M'$, the map $M'$ has $l$ connected components $M'_1,\ldots,M'_l$ and, for each $i$, the two extremities of $e_i$ belong to $M'_i$ and $M'_{i+1}$ (convention: $M'_{l+1}=M'_1$).\\
 
Choose any admissible loop ${L}$, it contains all the edges $e_i$. Indeed, if we look at a map of the kind $M'' = M \backslash E'$, with $E' \subsetneq \{e_1,\ldots,e_l\}$, all edges of $\{e_1,\ldots,e_l\} \backslash E'$ do not belong to any loop of $M''$ and are never erased in the computation of $D_1(M)$. So the only term in $T_{{L}}(M)$ which contribute to the coefficient of $M'$ is $(-1)^{l-1} M'$.
\end{proof}

\section{Decompositions and cumulants}\label{sect4}
In section \ref{sect3}, we have built an operator $D$ on bicolored labeled maps which leaves $N$ invariant and takes value in the ring spanned by forests. If we replace $N^{\dec{\tau}}$ by $N(D(M^{\dec{\tau}}))$ in the right hand side of equation (\ref{defkergen}), we obtain a decomposition of $\Sigma'_\mu$ as an algebraic sum of products of power series associated to minimal factorizations. In order to have something that looks like (\ref{existkergen}), we regroup some terms and make free cumulants appear through formula (\ref{stancum}). To do this, it will be useful to encode these associations of terms into combinatorial objects that we will call cumulant maps.\\

\comment{
In this section, we will construct thanks particular maps some functions $N_\varphi$ satisfying the condition (\ref{cumphi})
$$N_\varphi: \big\{(\dec{\tau}) \in S(k) \times S(k) \text{ s.t. } \tau \overline{\tau}=\sigma \big\} \longrightarrow \C[[\s{p},\s{q}]]$$
where the $c_i(\varphi)$ are cycles in $S(k)$ depending on $\varphi$. This map will be of the form $F \circ M_\varphi$. Then, we will associate to these functions integers $m_\varphi$ such that, for each pair $\tau,\bar{\tau} \in S(k)$,
$$D(M^{\dec{\tau}})=\sum_\varphi (-1)^{t_\varphi-1} m_\varphi M_\varphi({\dec{\tau}}).$$
So equation (\ref{Fdec}) will just be a corollary of the D-invariance of $N$
}

\subsection{Cumulant maps}
\begin{figure}
\begin{center}
\includegraphics{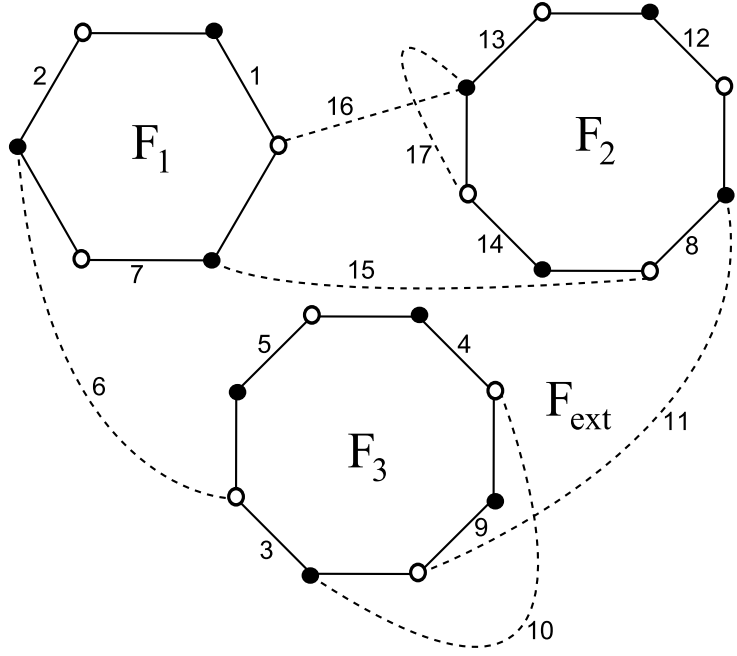}\caption{\footnotesize Example of a well-labeled cumulant map of resultant (1 \ldots 17)}
\label{figinjmap}
\end{center}
\end{figure}

\begin{Def}
A cumulant map $\M$ of size $k$ is a triple $(M_\M,\s{F},\iota)$ where $M_\M$ is a bicolored map with $|E|-|V|=k$, $\s{F}=(F_1,\ldots,F_t)$ is a family of faces of $M_\M$ such that
\begin{itemize}
\item The faces $F_1,\ldots,F_t$ are polygons (see definition \ref{defmaps})
\item Every vertex of $M_\M$ belongs to exactly one face among $F_1,\ldots,F_t$ ;
\end{itemize}
and $\iota$ is a function from $E \backslash \bigcup_i (E(F_i))$ (the set $E(F)$ was introduced in definition \ref{defmaps}) to $\N^\star$ (see figure \ref{figinjmap} for an example). As in the case of classical maps, if $\iota$ is a bijection of image $[k]$, the cumulant map is said well-labeled.\\

By definition, the number of connected components of $\M$ is the one of $M_\M$ and its resultant $\sigma_\M$ is the product of the cycles associated to the faces of $M_\M$ different from $F_1,\ldots,F_t$.
\end{Def}

\subsubsection{Non-crossing partitions as compressions of a polygon}
Consider a polygon with $2j$ vertices, alternatively black and white. We choose an orientation, begin at a black vertex and label the edges $1',1,2',2,\ldots,j',j$. Given a non-crossing partition $\pi \in NC(j)$, we glue, for each $i$, the edge $i$ with the edge $\sigma_\pi(i)'$ ($\sigma_\pi$ is the permutation of $[Id_j;(1 \ldots j)]$ canonically associated to $\pi$ by proposition \ref{bijpermncp}) so that their black extremities are glued together and also their white ones. In each of these gluings we only keep the label without $'$. The result is the labeled bicolored planar tree associated to the pair $\big(\sigma_\pi,\sigma_\pi^{-1}(1 \ldots j)\big)$.\\

This construction defines a bijection between $NC(j)$ and the different ways to compress a polygon with $2j$ vertices (with labeled edges) in a bicolored labeled planar tree with $j$ edges. So we reformulate (\ref{stancum}):
\begin{equation}\label{stancumarb}
R_{j+1}=\sum_{\substack{T \text{ tree obtained by compression}\\ \text{of a polygon of $2j$ vertices}}} (-1)^{|V_w(T)|+1} N(T),
\end{equation}
as power series in $\s{p}$ and $\s{q}$ (where $|V_w(T)|$ is the number of white vertices of $T$). If we consider a polygon without the labels $1',1,\ldots,j',j$, the bijection between $NC(j)$ and the different ways to compress it as a tree is only defined up to a rotation of the polygon but this formula is still true.\\

Given a cumulant map $\M$, consider all maps $M$ obtained from $M_\M$ by compressing each $F_i$ into a tree (we do not touch the edges - dotted in our example - which do not belong to any face $F_i$). Such maps $M$ have the same number of connected components as $\M$ and are maps of pairs of permutations whose product is the resultant of $\M$. The disjoint union of the trees obtained by compression of the face $F_i$ is a covering forest of $M$ with no trivial trees (\text{i.e.} with only one vertex), which is denoted $F_M$.\\

\begin{figure}
\begin{center}
\includegraphics{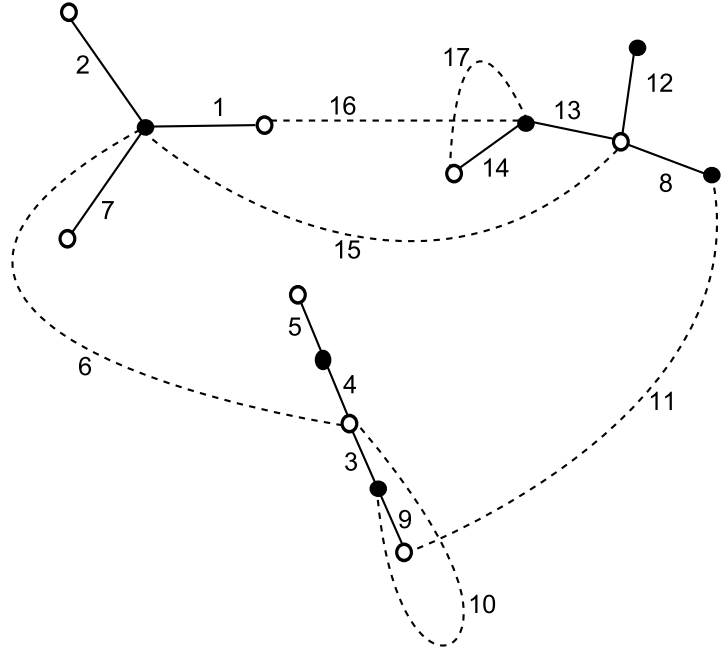}\caption{\footnotesize Example of a map obtained by compressing the polygons of the cumulant map of figure \ref{figinjmap}.}
\label{figcartintex}
\end{center}
\end{figure}

\begin{example}
The map $M$ of the figure \ref{figcartintex} can be obtained from the cumulant map of the figure \ref{figinjmap} by compressing each polygon into a tree in a certain way. The corresponding forest $F_M$ can be seen on the figure by erasing the dotted edges.\\
\end{example}

Let $\M$ be a cumulant map of resultant $\sigma$. Consider the function
$$N_\M: \big\{(\dec{\tau}) \in S(k) \times S(k) \text{ s.t. } \tau \overline{\tau}=\sigma \big\} \rightarrow \C[[\s{p},\s{q}]],$$
defined  by:
\begin{itemize}
\item If the map $M^{\dec{\tau}}$ is obtained from $M_\M$ by compressing in a certain way (necessarily unique) the faces $F_1,\ldots,F_t$, we put:
$$N_\M(\dec{\tau})=N\big(F_{M^{\dec{\tau}}}\big).$$
\item Else $N_\M(\dec{\tau})=0$.
\end{itemize}

This function fulfills:
\begin{equation}\label{cumphi}
\sum_{\substack{\dec{\tau} \in S(k)\\ \tau \overline{\tau} = \sigma_\M}} (-1)^{|C(\tau)|+t_\M}  N_\M({\dec{\tau}}) = \prod_{i=1}^{t_\M} R_{j_i+1}.
\end{equation}
\begin{proof}Use formula (\ref{stancumarb}) in the right hand side and expand it: the non-zero terms of the two sides of equality are exactly the same (with same signs because and $M$ and $F_M$ always have the same number of white vertices).
\end{proof}

Thanks to this property, this type of functions are a good tool to put series associated to forests together to make product of free cumulants appear.

\begin{Rem}\label{remint}
Let $\M$ be a cumulant map of resultant $\sigma$. The sets
$$\begin{array}{rc}
& \big\{\tau \in S(k) \text{ such that }N_\M(\tau,\tau^{-1}\sigma) \neq 0 \big\}\\
\text{and} & \big\{\overline{\tau} \in S(k) \text{ such that }N_\M(\sigma\overline{\tau}^{-1},\overline{\tau}) \neq 0 \big\}
\end{array}$$
are intervals $I_\M$ and $\overline{I_\M}$ of the symmetric group. So they are isomorphic as posets to products of non-crossing partition sets (for the order described in paragraph \ref{subsectncp}). The power series $N_\M(\tau,\tau^{-1}\sigma)$ is simply the one associated to the image of $\tau$ by this isomorphism (this image is defined up to the action of the full cycle on non-crossing partitions, so the associated power series is well-defined) and equation (\ref{stancum}) is a consequence of this fact.
\end{Rem}

%In the other sense, if we have a map $M_1$ of a pair $\dec{\tau}$ and a covering forest $M_0$ with no trivial trees (only one vertex), we can obtain a cumulant map $\M_{M_1,M_0}$ by replacing every tree $T_i$ by a polygon $F_i$ of length $2 |E(T_i)|$. This is the only cumulant map $\M$ for which:
%$$M_\varphi ({\dec{\tau}})\big)$$

%These applications have some nice properties:
%\begin{itemize}
%\item The number of connected components of $M_1$ is the one of $\M$ although the one of $M_2$ is the number $t_\M$ of polygons.
%\item The product $\tau \overline{\tau}$ is the product of $\M$.
%\end{itemize}

\subsection{Multiplicities}\label{cummapdec}
As for classical maps in paragraph \ref{comdec}, we define a decomposition operator for cumulant maps. Denote by $\mathbb{A}_{cm,\leq 1}$ the ring generated as $\Z$-module by the cumulant maps with at most one external half-edge by connected component. If $\M$ is a cumulant map, denote by $M'_\M$ the map obtained by repla\-cing, for each $i$, the face $F_i$ by a vertex (this map is not bicolored but each edge has a privileged orientation: the former white to black orientation).

\begin{Defprop}
There exists a unique linear operator $$\Decr: \mathbb{A}_{cm,\leq 1} \rightarrow \mathbb{A}_{cm,\leq 1}$$ such that:
\begin{itemize}
\item If $M'_\M$ has only one vertex, then $\Decr(\M)=\M$ ;
\item If $\M$ has more than one connected components ($\M=\prod_i \M^i$), then one has $\Decr(\M)=\prod \Decr(\M^i)$ ;
\item If $\M$ has only one connected component and no external half-edge, let $h$ be the half-edge of black extremity of its edge with the smallest label. We denote by $\overline{\M}$ the cumulant map obtained by adding one external half-edge between $h$ and its successor (as some edges have no labels, the half-edge is never in one of the faces $F_i$). Then $\Decr(\M)=\Decr(\overline{\M})$
\item If $M'_\M$ has only one connected component with one half-edge but no admissible loops, denote by $e_1,\ldots,e_l$ the edges leaving the same face $F_{i_0}$ as the external half-edge. The map $M_\M \backslash F_{i_0}$ has $l$ connected components $M_1,\ldots,M_l$, each with an external half-edge (at the place where $e_i$ leaves $M_i$). These maps have a cumulant map structure $M_i=M_{\M_i}$. Then $\Decr(\M)$ is given by the formula:
$$\Decr(\M)=\phi_\M(\Decr(\M_1),\ldots,\Decr(\M_l)),$$
where $\phi_\M$ is the multilinear operator on algebraic sums of sub-cumulant maps of the $\M_i$'s defined as $\phi_M$ in paragraph \ref{comdec}.
\item Else, consider $D_1(M'_\M)$ thanks to remark \ref{decnonbic}. In each map of the result, replace the vertices by faces $F_i$ and denote the resulting sum of cumulant map by $CM(D_1(M'_\M))$. Then,
$$\Decr(M)=\Decr(CM(D_1(M'_\M))).$$
\end{itemize}
\end{Defprop}

\begin{Def}
The multiplicity $c(\M)$ of a cumulant map $\M$ is the coefficient of the disjoint union of the faces $F_i$ in the decomposition $\Decr(\M)$ multiplied by $(-1)^{t_{\M}-1}$ (it can be zero!).
\end{Def}

Proposition \ref{propsign} is also true for cumulant maps and $\Decr$. So $c(\M)$ is non-negative if $\M$ is connected.\\

If $M$ is a map and $F_M$ a covering forest without trivial trees of $M$, denote by $\M_{M,F_M}$ the cumulant map obtained by replacing in $M$ each tree of $F_M$ by a polygon. The corresponding map $M'_{M,F_M}$ is
obtained from $M$ by replacing all trees of $F_M$ by a vertex. So the edges of $M \backslash F_M$ are in bijection with those of $M'_{M,F_M}$.

\begin{Lem}\label{lemcoef}
For any bicolored labeled map $M$, one has
$$D(M)=\sum_{F_M \subset M} (-1)^{t_{F_M}-1} c(\M_{M,F_M}) F_M,$$
where the sum runs over covering forests of $M$ with no trivial trees.
\end{Lem}

\begin{proof}
Let $F_M \subset M$ be a covering forest with no trivial trees of a bicolored labeled map. The operator $D$ applied to $M$ consists in making transformations of type $T_{{L}}$ with restricted choices until there are only forests left.
Thanks to remark \ref{ordre}, we choose loops containing a vertex of $T_\star$ (the tree of $F_M$ containing the external half-edge) as long as possible. As we are interested in the coefficient of $F_M$, we can forget at each step all maps that do not contain $F_M$. Now we notice that doing an elementary transformation with respect to $L$ and keeping only maps containing $F_M$ is equivalent to applying formula (\ref{transel}) with $E(L) \cap (M \backslash F_M)$ instead of $E(L)$.\\

%Thanks remark \ref{ordre}, we can in the induction definition of $D$ begin by making transformations with respect to loops containing a vertex of $T_\star$ (the tree of $F_M$ containing the half-edge).

As edges of $M \backslash F_M$ are in bijection with edges of $M'_{M,F_M}$, this new set of erasable edges is a set of edges of $M'_{M,F_M}$. With our choice of order of loops, this set of edges of $M'_{M,F_M}$ is always the set of erasable edges of an admissible transformation. So, computing $D(F_M)$ and keep only the submap containing $F_M$ is the same thing as computing $\Decr(\M_{M,F_M})$, except that we have trees instead of the polygonal faces. This shows that the coefficient of $F_M$ in $D(M)$ is the same as the one of the unions of the faces $F_i$ in $\Decr(\M_{M,F_M})$. The lemma is now obvious with the definition of the multiplicity of cumulant maps.
\end{proof}

With the notation of the previous paragraph, the lemma implies:
\begin{equation}\label{Fsumphi}
N(D(M^{\dec{\tau}}))=\sum_{\substack{\M \text{ cumulant maps}\\ \text{of resultant } \sigma}} (-1)^{t_\M-1} c(\M) N_\M(\dec{\tau}).
\end{equation}

\begin{Rem}
By remark \ref{remint} and lemma \ref{lemcoef}, for every $\sigma \in S(k)$, the family of intervals $I_\M$, where $\M$ describes the set of cumulant maps of resultant $\sigma$ with multiplicities $(-1)^{t_\M-1} c(\M)$, is a signed covering (the sum of multiplicities of intervals containing a given permutation is $1$) of the symmetric group by intervals $[\pi,\pi']$ such that
\begin{itemize}
\item The quantity $|C(\tau)|+|C(\tau^{-1}\sigma)|$ is constant on these intervals ;
\item The intervals are centered: $|C(\pi^{-1}\sigma)|=|C(\pi')|$.
\end{itemize}
Note that the power series $N$ does not appear in this result but is central in our cons\-truction. This interpretation of Kerov's polynomials' coefficients was conjecturally suggested by P. Biane in \cite{Bi2}.
\end{Rem}

\subsection{End of the proof of main theorem}
We use the $D$-invariance of $N$ to write $\Sigma'_{\mu}$ as an algebraic sum of power series associated to minimal factorizations:
\begin{eqnarray*}
\Sigma'_{\mu} & = & \sum_{\substack{\dec{\tau} \in S(k)\\ \tau \overline{\tau} = \sigma\\ <\dec{\tau}> \text{ trans.}}} (-1)^{|C(\tau)|+1} N(D(M^{\dec{\tau}}));\\
& = & \sum_{\substack{\dec{\tau} \in S(k)\\ \tau \overline{\tau} = \sigma\\ <\dec{\tau}> \text{ trans.}}} (-1)^{|C(\tau)|+1} \left[ \sum_{\substack{\M \text{ cumulant maps}\\ \text{of resultant } \sigma}} (-1)^{t_\M-1} c(\M) N_\M(\dec{\tau})\right].
\end{eqnarray*}
The second equality is just equation (\ref{Fsumphi}). Now, we change the order of summation (note that transitive factorizations have connected maps, so appear only as compressions of connected cumulant maps) and use (\ref{cumphi}):
\begin{eqnarray}
\Sigma'_{\mu} & = &\sum_{\substack{\M \text{ connected}\\ \text{cumulant map of}\\ \text{resultant }\sigma}} c(\M) \left[ \sum_{\substack{\dec{\tau} \in S(k)\\ \tau \overline{\tau} = \sigma}} (-1)^{|C(\tau)|+t_\M} N_\M(\dec{\tau})\right]; \nonumber\\
&=&\sum_{\substack{\M \text{ connected}\\ \text{cumulant map of}\\ \text{resultant }\sigma}} c(\M) \left[ \prod_{i=1}^{t_\M} R_{j_i(\M)+1} \right].
\end{eqnarray}
This ends the proof of theorem \ref{thkergen} because:
\begin{itemize}
\item the multiplicity of a connected cumulant map is non negative ;
\item the monomials in the $R_i$'s are linearly independent as power series in $\s{p}$ and $\s{q}$.
\end{itemize}

\section{Computation of some particular coefficients}\label{sect5}
\subsection{How to compute coefficients?}
In the proof of the main theorem, we have observed that the coefficient of the monomial $\prod\limits_{i=1}^t R_{j_i+1}$ in $K'_{\mu}$ is the sum of $c(\M)$ over all connected cumulant maps $\M$ of resultant $\sigma$, with $t$ polygons of respective sizes $2j_1,\ldots,2j_t$.\\

But it is easier to look, instead of the connected cumulant map $\M$, at the map $M_0$ obtained from $M_\M$ by compressing each polygon in a tree with only one black vertex. Recall that, in this context, $F_M$ is the disjoint union of these trees. Thanks lemma \ref{lemcoef}, the coefficient of $F_M$ in $D(M)$ is, up to a sign, equal $c(\M)$. Note that each pair $(M,F_M)$, where $M$ is the map of a transitive decomposition of $\sigma$ and $F_M$ a covering forest whose trees have exactly one black vertex and at least a white one, can be obtained this way from one cumulant map $\M$.\\

This remark leads to the following proposition, which will be used for explicit computations in the next paragraphs:
\begin{Prop}\label{calccoef}
The coefficient of monomial $\prod\limits_{i=1}^t R_{j_i+1}$ in $K'_{\mu}$ is the coefficient of the disjoint union of $t$ trees with one black and respectively $j_1,\ldots,j_t$ white vertices in
$$(-1)^{t-1} \sum_{\substack{\dec{\tau} \in S(k)\\ \tau \overline{\tau}=\sigma, <\dec{\tau}> \text{trans.} \\ |C(\overline{\tau})|=t}} D(M^{\dec{\tau}}).$$
\end{Prop}
As remarked before for coefficients of monomials of low degree, all the coefficients can be computed by counting some statistics on permutations in $S(k)$ (which can be much smaller than the symmetric group whose character values we are looking for).

\subsection{Low degrees in R}
\subsubsection{Linear coefficients}
A direct consequence of proposition \ref{calccoef} is the (well-known) combinatorial interpretation of coefficients of linear monomials in $\s{R}$: the coefficient of $R_{l+1}$ in $K'_{\mu}$ (or equivalently in $K_\mu$) is the number of permutations $\tau \in S(k)$ with $l$ cycles whose complementary $\overline{\tau}=\tau^{-1} \sigma$ is a full cycle, that is to say exactly the number of factorizations of $\sigma$, whose map has exactly one black vertex and $l$ whites. Indeed, if $M$ is a map with one black vertex, it is connected and has only loops of length $2$. So transformations with respect to these loops just consist in erasing an edge and $D(M)$ is a tree with one black vertex and as many white vertices as $M$.

\subsubsection{Quadratic coefficients} 
We have to compute $D(M)$, where $M$ is a connected map with two black vertices. Denote $w_0,\ldots,w_u$ the white vertices of $M$ linked to both black vertices. The first step is the computation of $D_1(\tilde{M})$, where $\tilde{M}$ is $M$ with an external half-edge $h$ (see definition \ref{defD}).\\

We begin by transformations with respect to all loops of length $2$ going through the extremity $\star$ of $h$. So we suppose that every $w_i$ is linked by only one edge $e_i$ to $\star$, but there can be more than one edge between $w_i$ and the other black vertex $v$, so we denote by $\s{f_i}$ the family of these edges. Let $e_i=\{h_i,h'_i\}$, where the extremity of $h_i$ is $\star$. With a good choice of numbering for the $w_i$, the cyclic order at $\star$ induces the order $h,h_0,\ldots,h_u$.

\begin{figure}
\begin{center}
\includegraphics{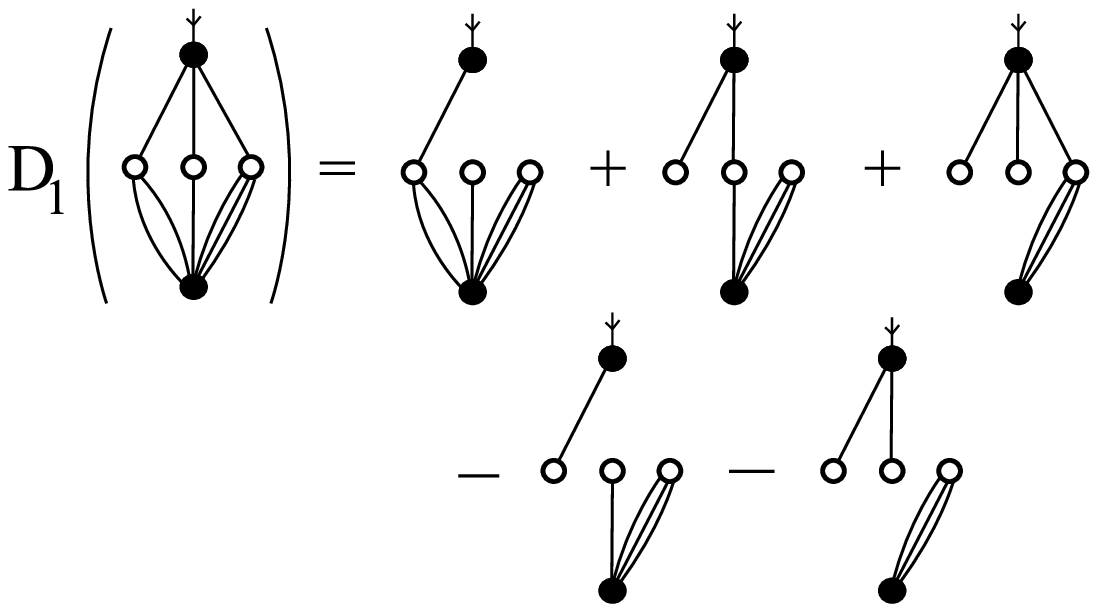}\caption{\footnotesize Elementary decomposition of a map with two black vertexes}
\label{figquad}
\end{center}
\end{figure}

\begin{Lem}
With these notations, we have:
\begin{equation}\label{dec2vb}
D_1(\tilde{M})=\sum_{i=0}^u \tilde{M} \backslash \{\s{f_0},\ldots,\s{f_{i-1}},e_{i+1},\ldots,e_u\} - \sum_{i=1}^u \tilde{M} \backslash \{\s{f_0},\ldots,\s{f_{i-1}},e_i,\ldots,e_u\}.
\end{equation}
An exemple for $u=3$ is drawn on figure \ref{figquad}.
\end{Lem}

\begin{proof}
If $u=0$, there is no admissible loop and this result is $D_1(\tilde{M})=\tilde{M}$. The case $u=1$ is left to the reader (it is an easy induction on the number of edge in $\s{f_0}$, the case where $\s{f_0}$ has two elements is contained in the case $H=2$ in the proof of definition-theorem \ref{invdecel}). Then we proceed by induction on $u$ by using the formula:
$$D_1(\tilde{M})=D_1 \big( D_1(\tilde{M} \backslash \{e_u\}) \cup \{e_u\} \big).$$
Suppose that lemma is true for $u-1$:
\begin{eqnarray}
D_1(\tilde{M} \backslash \{e_u\}) \cup \{e_u\}&=&\sum_{i=0}^{u-1} \tilde{M} \backslash \{\s{f_0},\ldots,\s{f_{i-1}},e_{i+1},\ldots,e_{u-1}\} \nonumber\\
\label{eqtechquad} & & \quad - \sum_{i=1}^{u-1} \tilde{M} \backslash \{\s{f_0},\ldots,\s{f_{i-1}},e_i,\ldots,e_{u-1}\}.
\end{eqnarray}
The graphs of the first line still have admissible loops. To compute their ima\-ge by $D_1$, we have to compute the image of the submaps whose set of edges is $\{e_i,\s{f_i},e_u,\s{f_u}\}$, since all other edges do not belong to any admissible loops. This is an application of the case $u=1$:
\begin{multline*}D_1(\tilde{M} \backslash \{\s{f_0},\ldots,\s{f_{i-1}},e_{i+1},\ldots,e_{u-1}\}) = \tilde{M} \backslash \{\s{f_0},\ldots,\s{f_{i-1}},\s{f_i},e_{i+1},\ldots,e_{u-1}\}\\ + \tilde{M} \backslash \{\s{f_0},\ldots,\s{f_{i-1}},e_{i+1},\ldots,e_{u-1},e_u\} - \tilde{M} \backslash \{\s{f_0},\ldots,\s{f_{i-1}},\s{f_i},e_{i+1},\ldots,e_{u-1},e_u\}.
\end{multline*}
Using this formula for each $i$, the first summand balances with the negative term in (\ref{eqtechquad}) (except for $i=u-1$) and the two other summands are exactly the ones in (\ref{dec2vb}). So the lemma is proved by induction.
\end{proof}

Now, in all maps appearing in $D_1(\tilde{M})$, there are only loops of length 2, so the end of the decomposition algorithm consists in erasing some edges without changing the number of connected components.\\

As explained in proposition \ref{calccoef}, we have to look at the sizes of trees in the two-tree forests (these forests come from the second sum of the right member of (\ref{dec2vb})). If, in $M$, there are $h^1_M$ white vertices linked to $\star$ (including the $w_i$) and $h^2_M$ to $v$, we obtain pairs of trees with $h^1$ and $h^2$ vertices, where $h^1$ and $h^2$ take all integer values satisfying the conditions:
$$\left\{\begin{array}{l} h^1-1 < h^1_M ; \\ h^2-1 < h^2_M ; \\ h^1+h^2=|V_w(M)|. \end{array} \right.$$
So any permutation with two black vertices contributes to coefficients of $R_{h^1} R_{h^2}$, where $h^1$ and $h^2$ verify the condition above. If $j \neq l$, a permutation may contribute twice to the coefficient of $R_j R_l$ if the conditions above are fulfilled for $j=h^1, l=h^2$ and for $l=h^1, j=h^2$. Finally, one has:
\begin{multline*}[R_j R_l] K_k = \left\{\begin{array}{cl} 1 & \text{if }j \neq l\\
1/2 & \text{if }j = l \end{array} \right\} \cdot \\
\qquad \sum_{\substack{\dec{\tau} \in S(k)\\ \tau \overline{\tau}=\sigma, <\dec{\tau}> \text{trans.} \\ |C(\overline{\tau})|=2}} \delta_{j \leq h^1_{M^{\dec{\tau}}}} \delta_{l \leq h^2_{M^{\dec{\tau}}}} + \delta_{l \leq h^1_{M^{\dec{\tau}}}} \delta_{j \leq h^2_{M^{\dec{\tau}}}},\end{multline*}
which is exactly the second part of theorem \ref{thlinquad} (the second $\delta$ in the equation above disappears if we consider permuations with numbered cycles).

\subsection{High degrees in p,q}
If the graded degree in $\s{p}$ and $\s{q}$ is high, the maps we are dealing with have few loops. Therefore, it is easier to compute their image by $D$ and to count them.

\begin{proof}[Proof of theorem \ref{thgenre1gen}]
Let $r,s,t,j_1,\ldots,j_t$ be integers such that $\sum j_i=r+s$. As in the whole paper $\sigma \in S(k)$ is a permutation of type $\mu$ (here $r,s$). We can suppose that $1$ is in the support of the cycle $c_1$ of $\sigma$ of size $s$.\\

We have to count connected maps with $r+s$ edges and $r+s$ vertices, that is to say, up to a change of orientation, one loop $L$. So, eventually by replacing $L$ by $L'$ (if $1$ is in the word associated to the external face, $L$ must be going counterclockwise), $D(M)=T_L(M)$. Only maps $M$ such that, in $D(M)$, there is (at least) a forest with one black vertex per tree, contribute to coefficients of Kerov's polynomials. In such maps, all vertices of $M \backslash {L}$ are white and only the forest $M \backslash E(L)$ (see formula (\ref{transel})) satisfies the condition above.\\

Let us consider such a map $M$. We can choose arbitrarily a first black vertex $b_1$ of $M$ ($M$ will be said marked) and number $b_1,\ldots,b_t$ all its black vertices in the order of $L$. Suppose that there are $w_i$ white vertices of $M \backslash{L}$ linked to $b_i$. Then $M$ contributes only to the coefficient of $\prod R_{w_i+2}$ in $K'_{f_1,f_2}$ (where $2f_1$ and $2f_2$ are the lengths of the two faces of $M$) with coefficient $1$.\\

We count the number of marked labeled maps $M$ contributing to the coefficient of $\prod\limits_{i=1}^t R_{j_i}$ in $K'_{r,s}$. They are of the form of the figure \ref{figgenre1r2} with:
\begin{itemize}
\item The word $(r_1 + s_1, r_2 + s_2,\ldots,r_t + s_t)$ is equal up to a permutation to $(j_1-2,\ldots,j_t-2)$
\item The length $r_1 + r_2 + \ldots +r_t$ of the face $F_r$ which is on the left side of ${L}$, is equal to $r$.
\end{itemize}

\begin{figure}
\begin{center}
\includegraphics{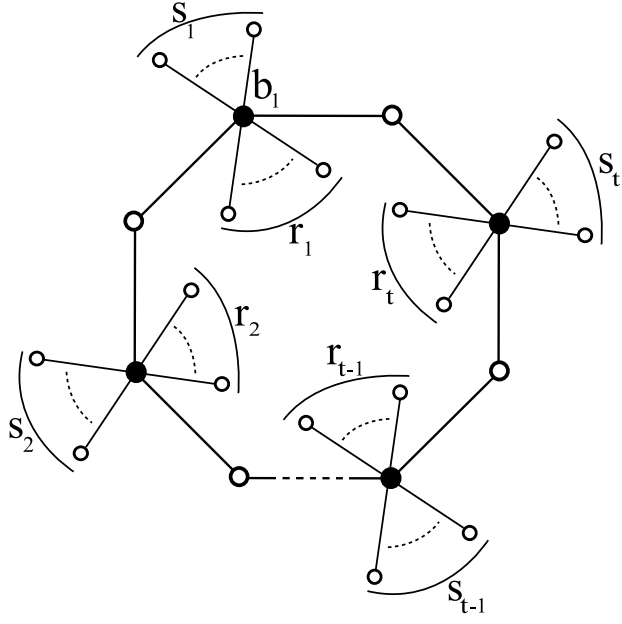}\caption{\footnotesize Maps contributing to terms of graded degree $r+s$ in $K'_{r,s}$.}
\label{figgenre1r2}
\end{center}
\end{figure}

Such a map can be labeled of $r \cdot s$ different ways such that its faces are the cycles of $\sigma$. Indeed, if we fix one element in the support of each cycle of $\sigma$, such a labeling is determined by the edges labeled by these elements. We have $r$ (resp. $s$) choices for the first (resp. second) one: the $r$ (resp. $s$) edges whose labels are in the word associated to the face $F_r$ (resp. $F_s$). As we deal for the moment with maps with a marked black vertex, all the numberings give a different map.\\

If we choose a permutation $\s{j'}-2$ of the word $(j_1-2,\ldots,j_t-2)$, non-negative integers $r_1,s_1,\ldots,r_t,s_t$ such that $\sum\limits_i r_i =r-t, \sum\limits_i s_i =s-t$ and $\s{r}+ \s{s}=\s{j'}-2$ and labels on the corresponding map, we obtain a marked map $M$ contributing to the coefficient of $\prod\limits_{i=1}^t R_{j_i}$ in $K'_{r,s}$. To obtain the number of such non-marked maps, we have to divide by $t$ (thanks to the labels, there is no problem of symmetry).\\

So the coefficient of $\prod\limits_{i=1}^t R_{j_i}$ in  $K'_{r,s}$ is
$$\frac{r \cdot s}{t} \ \Perm(\s{j}) \ \big|\big\{(r_1,s_1,\ldots,r_t,s_t)\big\}\big|,$$
where $r_1,s_1,\ldots,r_t,s_t$ describe the set of non-negative integers satisfying the equations
$$\left\{\begin{array}{c} 
r_1 + s_1 = j_1 - 2 ;\\
\vdots\\
r_t + s_t = j_t - 2 ;\\
r_1 + \ldots + r_t = r - t.
\end{array}\right.$$

But, in the system of equations satisfied by the $r_i$'s and the $s_i$'s, we can forget the $s_i$'s and only keep an inequality on each $r_i$ ($r_i \leq j_i -2$), which corresponds to the positivity of $s_i$. So the cardinal of the set in the formula above is exactly $N(j_1 - 2,\ldots,j_t - 2 ; r - t)$\\
\end{proof}

We use the same ideas for subdominant term in the case $l(\mu)=1$.

\begin{proof}[proof of theorem \ref{thgenre1}]
To compute the coefficients of a monomial of graded degree $k-1$ in $K_k$, we have to count the contributions of labeled maps with $k$ edges, $k-1$ vertices and one face. As in the previous proof, if a map has a non-zero contribution, all vertices which do not belong to any loop are white. Such maps can be sorted in five classes: see figure \ref{figgenre1r1} for types $a$ and $b$, type $c$ (resp. $d$) is type $b$ with one black and one white (resp. two white) vertices at the extremities and type $e$ is type $a$ with a white central vertex of valence 4 instead of a black one.\\
\begin{figure}
\begin{center}
\includegraphics{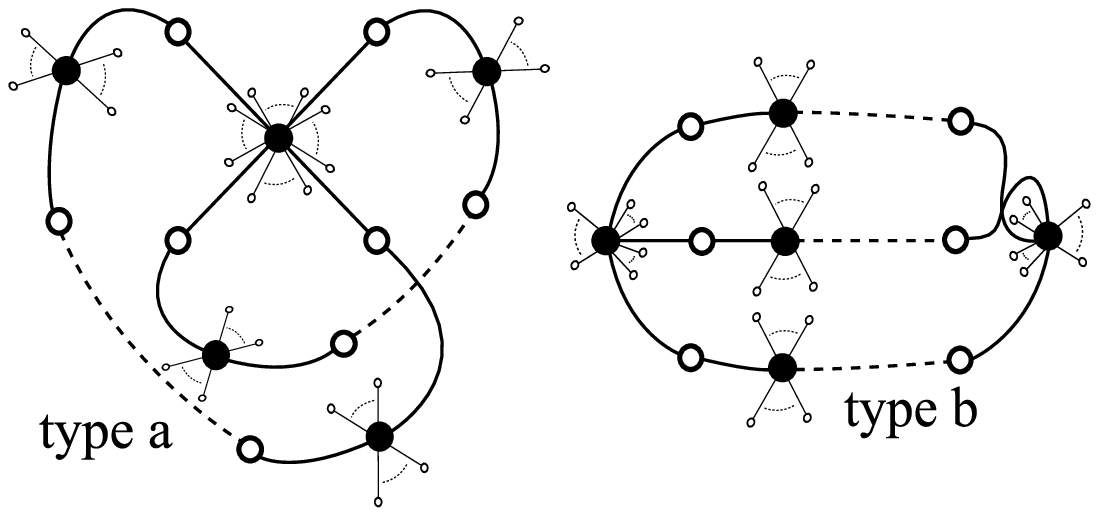}\caption{\footnotesize Maps contributing to terms of graded degree $k-1$ in $K_k$.}
\label{figgenre1r1}
\end{center}
\end{figure}

Thanks to the case $H=2$ in the proof of definition-theorem \ref{invdecel}, the decomposition of these maps is easy to compute:
\begin{description}
	\item[Types $a$ and $e$] The two loops have no edges in commun and their associated transformations commute ;
	\item[Types $b$, $c$ and $d$] We obtain a result close to the one of figure \ref{figexinv}.
\end{description}
Here is the description of the forests with $t$ trees for each type (it is quite surprising that it does not depend on the labels).
\begin{description}
\item[Type $a$] In $D(M)$, there is one forest $F$ with one black star per tree: in addition to those which do not belong to loops, there are two white vertices linked to the central black vertex and one to each other black vertex.
\item[Type $b$] In $D(M)$, there are two forests $F_1$ and $F_2$ with one black star per tree: in $F_1$ (resp. in $F_2$), in addition to those which do not belong to loops, there are two white vertices linked to the vertex at the left (resp. right) extremity and one to each other black vertex (including the right (resp. left) extremity).
\item[Type $c$] In $D(M)$, there is one forest $F$ with one black vertex per tree: in addition to those which do not belong to loops, there is one white vertex linked to each black vertex.
\item[Types $d$ and $e$] In $D(M)$, there is no forest $F$ with one black vertex per tree.
\end{description}

Now we compute the coefficient of $\prod\limits_{i=1}^t R_{j_i}$ in $K_k$. We give all the details only for the contributions of maps of type $a$.\\

If we mark an half-edge of extremity the central black vertex in a map $M$ of type $a$, we number the black vertices of $M$ by following the face of $M$ beginning by this half-edge (but not by the central black vertex). As in the previous proof, a map contributing to this monomial with a marked half-edge of extremity the central black vertex ($4$ choices) is given by:
\begin{itemize}
\item A permutation $\s{j'}$ of the word $(j_1,\ldots,j_t)$ ($j'_i$ is the number of vertices of the tree of $F$ of black vertex $b_i$).
\item The length of the first loop, \text{i.e.} the label $p \in [t]$ of the central black vertex. 
\item For each black vertex different from the central one, we have to link $j'_i-2$ white vertices that do not belong to loops. We have to fix the number of these vertices which are on a given side of the loop: there is $j'_i-1$ possibility.
\item Idem for the central black vertex except that we have $j_p'-3$ white vertices to place in $4$ sides, so $\left(\begin{array}{c} j_p'\\3 \end{array}\right)$ possibilities.
\item The labels of such a map are determined by the choice of one edge which has the label $1$, so $k$ possibilities.
\end{itemize}

Finally the contribution of type $a$ maps to the coefficient of $\prod\limits_{i=1}^t R_{j_i}$ in $K_k$ is
$$C_a = \frac{k}{4} \sum_{\s{j'}} \left[\sum_{p=1}^t \frac{j'_p(j'_p-2)}{6} \prod_{i=1}^t (j'_i - 1)\right].$$
The expression in the bracket is symmetric in $\s{j'}$, so equal to its value for $\s{j}$:
$$C_a = \frac{k}{4} {|\Perm(\s{j})|} \prod_{i=1}^t (j_i - 1) \sum_{p=1}^t \frac{j_p(j_p-2)}{6}.$$\\

We can find similar arguments for types $b$ and $c$:
\begin{itemize}
\item In type $b$, $p_1$ and $p_2$ are the labels of the black vertices at the extremities if we numbered by following the face beginning just after an extremity ($6$ possibilities to choose where to begin) ;
\item In type $c$, $p_1$ is the label of the black extremity and $p_2$ of the black vertex preceding the white extremity if we begin just after the white extremity ($3$ possibilities to choose where to begin), note also that in this type we have to symmetrize our expression in  $\s{j'}$.
\end{itemize}
We obtain:
\begin{eqnarray*}
C_b & = & \frac{k}{6} |\Perm(\s{j})| \prod_{i=1}^t (j_i - 1) \sum_{1 \leq p_1 < p_2 \leq t} \frac{j_{p_1}(j_{p_2}-2)}{4} +\frac{j_{p_2}(j_{p_1}-2)}{4};\\
C_c & = & \frac{k}{3} |\Perm(\s{j})| \prod_{i=1}^t (j_i - 1) \sum_{1 \leq p_1 \leq p_2 \leq t} \frac{1}{2} \left(\frac{j_{p_1}}{2} + \frac{j_{p_2}}{2}\right).
\end{eqnarray*}\\

Finally, if we note 
$$A = \frac{k}{24} |\Perm(\s{j})| \prod_{i=1}^t (j_i - 1),$$
and split the summation in $C_c$ into the cases $j_{p_1} < j_{p_2}$ and $j_{p_1}=j_{p_2}$, the coefficient we are looking for is:
\begin{eqnarray*}
C_a+C_b+C_c & = & A \left(\sum_{p=1}^t j_p(j_p-2) + \sum_{1 \leq p_1 \leq t} 4 j_{p_1}\right.\\
& &\qquad \left. + \sum_{1 \leq p_1 < p_2 \leq t} \big( j_{p_1}(j_{p_2}-2) + j_{p_2}(j_{p_1}-2) + 2j_{p_1} + 2j_{p_2} \big) \right) ;\\
& = & A \left(2 \sum_{p=1}^t j_p + \sum_{p=1}^t j_p^2 + \sum_{1 \leq p_1 < p_2 \leq t} \big( j_{p_1}j_{p_2} + j_{p_2}j_{p_1} \big)\right);\\
& = & A \left[ \left(\sum_{p=1}^t j_p \right)^2 + 2\sum_{p=1}^t j_p \right];\\
& = & A \big( (k-1)^2 + 2(k-1) \big) = A(k-1)(k+1),
\end{eqnarray*}
which is exactly the expression claimed in theorem \ref{thgenre1}.
\end{proof}

\section*{Acknowledgements}
The author would like to thank his adviser P. Biane for introducing him to the subject, helping him in his researches and reviewing (several times) this paper. He also thanks Piotr \'Sniady for stimulating discussions.

\end{document}